\documentclass[]{elsarticle}
\pdfoutput=1

\makeatletter
\def\ps@pprintTitle{
\let\@oddhead\@empty
\let\@evenhead\@empty
\def\@oddfoot{\centerline{Preprint accepted to the Elsevier Journal of Parallel Computing, Dec 2017.}}
\let\@evenfoot\@oddfoot}
\makeatother

\usepackage[margin=2.0cm]{geometry}
\usepackage{amsfonts}
\usepackage{amsmath}
\usepackage[dvipsnames]{xcolor}
\usepackage{amsfonts}
\usepackage{graphicx}
\usepackage{algorithmic}
\usepackage{amsmath}
\usepackage[justification=centerlast]{caption}
\usepackage{subcaption}
\usepackage{booktabs}
\usepackage{pgfplots}
\usepackage{pgfplotstable}
\pgfplotsset{compat=1.12}
\usepackage{multirow}
\usepackage{colortbl}
\usepackage{mathabx}
\usepackage{algorithm,caption}
\usepackage{listings}
\usepackage{float}

\begin{document}

\begin{frontmatter}
	\title{Accelerated Cyclic Reduction: A Distributed-Memory \\ Fast Solver for Structured Linear Systems}
	\author[kaust-ecrc]{Gustavo Ch{\'a}vez\corref{joris}}
	\ead{gustavo.chavezchavez@kaust.edu.sa}
	\author[aub-cs]{George Turkiyyah}
	\author[kaust-ecrc]{Stefano Zampini}
	\author[kaust-ecrc]{Hatem Ltaief}
	\author[kaust-ecrc]{David Keyes}

	\cortext[joris]{Corresponding author}
	\address[kaust-ecrc]{King Abdullah University of Science and Technology (KAUST).}
	\address[aub-cs]{Department of Computer Science, American University of Beirut (AUB).}

	\begin{abstract}
	We present Accelerated Cyclic Reduction (ACR), a distributed-memory fast solver for rank-compressible block tridiagonal linear systems arising from the discretization of elliptic operators, developed here for three dimensions.
	Algorithmic synergies between Cyclic Reduction and hierarchical matrix arithmetic operations result in a solver that has $O(k~N \log N~(\log N + k^2))$ arithmetic complexity and $O(k~N \log N)$ memory footprint, where $N$ is the number of degrees of freedom and $k$ is the rank of a block in the hierarchical approximation, and which exhibits substantial concurrency.
	We provide a baseline for performance and applicability by comparing with the multifrontal method with and without hierarchical semi-separable matrices, with algebraic multigrid and with the classic cyclic reduction method.
	Over a set of large-scale elliptic systems with features of nonsymmetry and indefiniteness, the robustness of the direct solvers extends beyond that of the multigrid solver, and relative to the multifrontal approach ACR has lower or comparable execution time and size of the factors, with substantially lower numerical ranks.
	ACR exhibits good strong and weak scaling in a distributed context and, as with any direct solver, is advantageous for problems that require the solution of multiple right-hand sides.
	Numerical experiments show that the rank $k$ patterns are of $O(1)$ for the Poisson equation and of $O(n)$ for the indefinite Helmholtz equation.
	The solver is ideal in situations where low-accuracy solutions are sufficient, or otherwise as a preconditioner within an iterative method.
	\end{abstract}

\end{frontmatter}

\section{Introduction}
\label{sec:introduction}

Cyclic reduction, introduced in \cite{hockney65}, is a direct solver for tridiagonal linear systems. It is effective for the solution of (block) Toeplitz and (block) tridiagonal matrices that arise from the discretization of elliptic PDEs \cite{buzbee70,gander97}. For the constant-coefficient Poisson equation, since each of the blocks of the discretized system is Fourier diagonalizable, cyclic reduction can be used in combination with the fast Fourier transform (FFT) to deliver optimal complexity, as proposed in the FACR method \cite{swarztrauber77}. However, in the presence of variable coefficients, the FFT-enabled version of cyclic reduction can not be used. The purpose of this work is to address the time and memory complexity growth in the presence of heterogeneous blocks with a variant called Accelerated Cyclic Reduction (ACR). The main observation is that elliptic operators have a hierarchical structure of off-diagonal blocks that can be approximated with low-rank matrices. Thus we approximate appropriate blocks of the initially sparse matrix with hierarchical matrices and operate on these blocks with hierarchical matrix arithmetics, instead of the usual dense operations, to obtain a direct solver of log-linear arithmetic and memory complexities. This philosophy follows recent work discussed below, but to our knowledge, this is the first demonstration of the utility of complexity-reducing hierarchical substitution in the context of cyclic reduction.

Cyclic reduction can be thought of as a direct Gaussian elimination on a permuted system that recursively computes the Schur complement of half of the unknowns until a single block remains or the system is small enough to be inverted directly. Schur complement computations have a complexity that is dominated by the cost of the inverse; by applying a red/black re-ordering of the unknowns, the linear system separates into two halves with block diagonal structure.
This decoupling addresses the most expensive step of the Schur complement computation regarding operation complexity and does so in a way that launches independent subproblems. This concurrency feature, in the form of recursive bisection, can be naturally implemented in a distributed-memory parallel environment. The stability of the block cyclic reduction has been studied in \cite{yalamov1996stability}, where the author presents error bounds for strictly and nonstrictly diagonally dominant matrices.

In order to simplify the description of the algorithm, in this work we consider structured linear systems arising from the discretizations or scalar PDEs on three-dimensional Cartesian grids. For three-dimensional problems of size $N=n^3$, where $n$ is the number of discretization points in the linear dimension of the target domain, the synergy of cyclic reduction and hierarchical matrices leads to a parallel fast direct solver of $O(k~N \log N~(\log N + k^2))$ arithmetic complexity, and $O(k~N \log N)$ memory footprint, where $k \ll N$ represents the numerical rank of compressed blocks. This is in contrast to $O(N^2)$ and $O(N^{1.5})$ respectively, if hierarchically low-rank matrices matrices are not used.

In this manuscript, we present ACR and its distributed-memory implementation, and we demonstrate its performance on a set of problems with various symmetry and spectral properties in three dimensions. These problems include the Poisson equation, the convection-diffusion equation, and the indefinite Helmholtz equation. We show that ACR is competitive in memory consumption and time-to-solution when compared to methods that rely on a global factorization and do not exploit the cyclic reduction structure.

\subsection{Related work}

Recent years have seen increasing interest in the use of hierarchical low-rank approximations to accelerate the direct solution of linear systems. In this section, we briefly describe some of this literature focusing primarily on efforts that target distributed-memory environments.

Arguably the most common approach for using hierarchical matrix representations in matrix factorizations is to use low-rank approximations to compress the dense frontal blocks that arise in the multifrontal variant of Gaussian elimination. The enabling property is that under proper ordering, many of the off-diagonal blocks of the Schur complement of discretized elliptic PDEs have an effective low-rank approximation \cite{chandrasekaran2010} that improves the memory and arithmetic estimates of conventional multifrontal solvers \cite{Duff83}. Furthermore, there are efficient low-rank approximation methods to perform the necessary arithmetic operations and preserve the low-rank representation during the factorization and solution stages of the solver. Within this general approach, various methods that differ in the particular data-sparse format used and in the algorithms for the computation of low-rank approximations have been developed.

In Wang et al. \cite{Wang2016} the authors investigate the use of the HSS format \cite{Vandebril05} to accelerate the parallel geometric multifrontal method, which results in a method known as the HSS-structured multifrontal solver (HSSMF). The general approach uses intra-node parallel HSS operations within a distributed-memory implementation of the multifrontal sparse factorization. This approach lowers the complexity of both arithmetic operations and memory consumption of the resulting HSS-structured multifrontal solver by leveraging the underlying numerically low-rank structure of the intermediate dense matrices appearing within the factorization process, driven by an optimal nested dissection ordering.

In a similar line of work, Ghysels et al. \cite{Ghysels15} also investigate a combination of the multifrontal method and the HSS-structured hierarchical format, extending the range of applicability of the solver to general non-symmetric matrices. Using the task-based parallelism paradigm, they introduce randomized sampling compression \cite{martinssonRandomized2011} and fast ULV HSS factorization \cite{xia_ulv_2012}. Under the assumption of the existence of an underlying low-rank structure of the frontal matrices, randomized methods deliver almost linear complexity; this reduces the asymptotic complexity of the solver, which is mainly attributed to the frontal matrices near the root of the elimination tree. The effectiveness of these task-based algorithms in combination with a distributed-memory implementation of the multifrontal method is available in an early stage software release of the package STRUMPACK \cite{Ghysels15}, which we will consider in the numerical experiments section of this article. The HSS format assumes a weak admissibility condition (see section \ref{sec:H_matrix_construction}), which in practice requires the use of large numerical ranks even for approximations with modest relative accuracy. Consequently, this stresses the memory requirements and increases overall execution time.

The hierarchical interpolative factorization \cite{ho2015_hif_de, ho2015_hif_ie} is another method for finding low-rank approximations that has proved to be a fast solver for symmetric elliptic PDEs and integral equations. This decomposition relies on a ``skeletonization'' procedure to eliminate a redundant set of points from a symmetric matrix to further compress the dense fronts. The key step in skeletonization uses the interpolative decomposition of low-rank matrices to achieve a quasi-linear overall complexity in factorization. The performance of hierarchical interpolative decomposition in a distributed-memory environment is reported in \cite{li2016distributed}.

A fast direct method for high-order discretizations of elliptic PDEs has been proposed by Martinsson et al. \cite{Martinsson2013460,gillman2014direct,Hao2016419}. The method is based on a multidomain spectral collocation discretization scheme and a hierarchy of nested grids, similar to nested dissection. It exploits analytical properties of elliptic PDEs to build Dirichlet-to-Neumann operators, by hierarchically merging these operators originating from smaller grids. When computations are done using the HSS data-sparse format, an asymptotic complexity of $O(N^{4/3})$ can be reached. The high-order discretizations used in this method makes it quite powerful in practice as they allow it reach the same accuracy with fewer degrees of freedom compared to second order discretizations. A distributed-memory implementation of this algorithm is in progress.

Even though this approach has larger asymptotic estimates than the log-linear performance of the methods above, because of the high-order discretization of the PDE, this method is quite powerful in practice as they can reach the same accuracy with fewer degrees of freedom as compared to second order discretizations. A distributed-memory implementation of this algorithm is in progress.

The BLR format \cite{weisbecker2013improving} has also been used to compress blocks into low-rank approximations to accelerate the factorization process of the multifrontal method. This format is compatible with numerical pivoting and is well-suited for the reuse of existing high-performance implementations of dense linear algebra kernels. Even though this format is not hierarchical, it has proven to be useful for a wide range of problems \cite{amestoy2015improving} within the distributed-memory implementation of the multifrontal method provided by the MUMPS library \cite{MUMPS1,MUMPS2}.

Rather than compressing and identifying individual blocks of the decomposition, another hierarchy-exploiting approach considers the system as a whole and seeks to construct a holistic decomposition of the full linear system. An example of such decomposition is the recursive computation of the inverse of a hierarchical matrix \cite{Izadi2012,Ambikasaran2013}, or the computation of its Cholesky or LU factorization \cite{ibra07,blackBoxHLU08}. These methods have generally much higher prefactors than methods that compress individual matrix blocks of the factorizations and are not usually competitive for large-scale problems; as an example, we refer the reader to \cite{Izadi2012} for a discussion of the challenges of scaling the construction of the inverse of a hierarchical matrix.

Pouransari et al. approximate fill-in via low-rank approximations with the $\mathcal{H}^2$ format; see \cite{Pouransari2016}. This format guarantees linear complexity provided that blocks correspond to well-separated clusters and have a data-sparse property. The algorithm starts by recursively bisecting the computational domain, implicitly forming a binary tree. The leaf nodes correspond to independent subdomains, and the internal nodes correspond to Schur complements to computed with low-rank arithmetic operations. The bottom-up elimination process is performed with a procedure referred to as ``extended sparsification'' in which the original matrix dimension grows by introducing auxiliary variables but nonetheless remains sparse. Alternatively, elimination can be performed with an in-place algorithm that keeps the matrix size constant. A related method with similar strategies as in this work is the so-called ``compress and eliminate'' solver \cite{Sushnikova2016}. A recent extension of this line of work into a distributed memory environment documented in \cite{Chen17}, demonstrates that concurrent processors can work on independent subdomains defined by their corresponding subgraphs, where interior vertices are eliminated concurrently. Communication is needed at the boundary vertices, but additional concurrency at the boundary is exploited trough graph coloring.

\subsection{Contributions}

The contribution of this work is the development of a parallel, robust and efficient method for the solution of block tridiagonal linear systems, with emphasis on systems that arise from the discretization of elliptic PDEs. ACR is a fast solver in the sense that it has a log-linear arithmetic complexity in operations count and memory consumption. The algorithm arrives at the solution in a finite number of steps, rather than iteratively converging to a solution, which makes it a direct solver with a tunable accuracy. The fact that ACR is entirely algebraic extends its range of applicability to problems with an arbitrary coefficient structure including nonsymmetry within the block tridiagonal sparsity structure, subject to their amenability to rank compression.  This entirely algebraic property gives the method robustness on problems that are challenging for iterative methods, while still maintaining asymptotic efficiency.

Two key features of the algorithm from a computational perspective are the simplicity of its parallelization and the regularity of its communication patterns in a distributed memory environment. The communication pattern is well-established beforehand and it is based on recursive bisection, as opposed to nested dissection with different block sizes at different levels of the factorizations. The amount of inter-node concurrency is proportional to the size of the blocks and it fits readily into a distributed-memory parallel environment. The algorithm also exhibits substantial intra-node concurrency, both in processing multiple blocks and within its hierarchical operations on individual blocks, which fits the multi-core architecture of modern supercomputers.

We demonstrate that our implementation is well suited for modern parallel multi-core systems and scalable in a distributed-memory environment. We also compare our implementation against other state-of-the-art direct solvers over a relevant class of problems and show competitive time to solution and memory requirements.

\section{Preliminaries} 
In this section, we review the building blocks of the proposed solver, namely hierarchical low-rank approximations and the cyclic reduction algorithm. 

\subsection{Hierarchical matrices}

A hierarchical matrix is a data-sparse representation that enables fast linear algebraic operations by using a hierarchy of off-diagonal blocks, each represented by a low-rank approximation, that can be tuned to guarantee an arbitrary precision. The approximation, sometimes referred to as compression, is performed via singular value decomposition, or with a related method that delivers a low-rank approximation with fewer arithmetic operations than the traditional SVD method. For the representation to be effective in terms of arithmetic operations and memory requirements, numerical ranks significantly smaller than the sizes of the various matrix blocks are required.

There are several hierarchical and non-hierarchical low-rank approximation formats available in the literature. In this work, we consider the $\mathcal{H}$-matrix format introduced by Hackbusch et al. in \cite{hackbusch99}. Being modular by design, ACR is not limited to the $\mathcal{H}$-format. In fact, the use of the $\mathcal{H}^2$-format would immediately translate to an additional reduction of one logarithmic factor in terms of arithmetic and memory complexity estimates, from $O(k~N \log N~(\log N + k^2))$ to $O(k~N \log N)$ in terms of operations, and $O(k^2~N \log N)$ to $O(k~N)$ in terms of memory requirements, however, we require a complete set of hierarchical matrix operations and fast construction, which at the time of this publication is still ongoing work within our group. Our implementation uses the $\mathcal{H}$-format arithmetic and its arithmetics operations provided by the HLibPro library. We refer to the reader to \cite{kriem05,Kriemann2014} for a discussion of the shared-memory scalability of these hierarchical matrix operations, their relative costs, and their performance on modern manycore architectures. HLibPro does not feature a distributed memory solver. We use its shared-memory kernels in combination with MPI to orchestrate parallel workload across nodes in a distributed memory environment as we will discuss in section \ref{parallelACR}.

\subsubsection{$\mathcal{H}$-matrix construction}
\label{sec:H_matrix_construction}

The structure of a hierarchical matrix in the $\mathcal{H}$ format can be described by four components: an index set, a cluster tree, a block cluster tree, and the choice of an admissibility condition. The index set $\mathcal{I} = \{0,1,\dots,N-1\}$ represents the number of degrees of freedom $N$. The cluster tree represents row/column groupings, and it is constructed by recursively subdividing the index set.  Once the cluster tree is formed, the block cluster tree defines matrix sub-blocks over the index $\mathcal{I} \times \mathcal{I}$. Its leaves are either low-rank blocks or small dense ones.  Finally, the admissibility condition determines whether a given block should be represented as a low-rank approximation or a dense block\footnote{The word ``block'' is overloaded in this discussion. It is used to denote the partitions of the block tridiagonal coefficient matrix of the problem. It is also used to denote the partitioning of a matrix into low-rank and dense subdivisions. When necessary to avoid confusion, we will use the word ``plane'' or ``plane-block'' to refer to the first meaning.}.

The first step for the construction of an $\mathcal{H}$-matrix is the definition of the cluster tree of unknowns. In this work, since each block row of the sparse matrix represents a plane from a three-dimensional regular discretization, we leverage the geometry information by selecting a binary space partitioning strategy to cluster the unknowns considering the two-dimensional domain representing the planes. 

The next step is the definition of a block cluster tree for these two-dimensional domains, which together with the admissibility condition determines the structure of the hierarchical representation of the plane-block. We chose a standard admissibility condition, as opposed to the simpler weak admissibility condition, because it provides the flexibility of selecting a range of coarser or finer blocks, tuned by an admissibility parameter $\eta$. Weak admissibility refers to a matrix decomposition where the $(1,2)$ and $(2,1)$ blocks are single low-rank blocks and the $(1,1)$ and $(2,2)$ blocks are recursively decomposed in a similar way. 
On the other hand, standard admissibility allows a more refined blocking of the matrix; the $\eta$ parameter appears in the inequality $min( diameter(\tau), diameter(\sigma) ) \leq \eta \cdot distance(\tau,\sigma)$, where $\tau$ and $\sigma$ denote two geometric regions defined as the convex hulls of two separate point sets $t$ and $s$ (nodes in cluster tree). A matrix block $A_{ts}$ satisfying the previous inequality is represented in a low rank form.

\begin{figure}[H]
	\begin{subfigure}{.5\textwidth}
		\centering
		\includegraphics[width=0.9\linewidth]{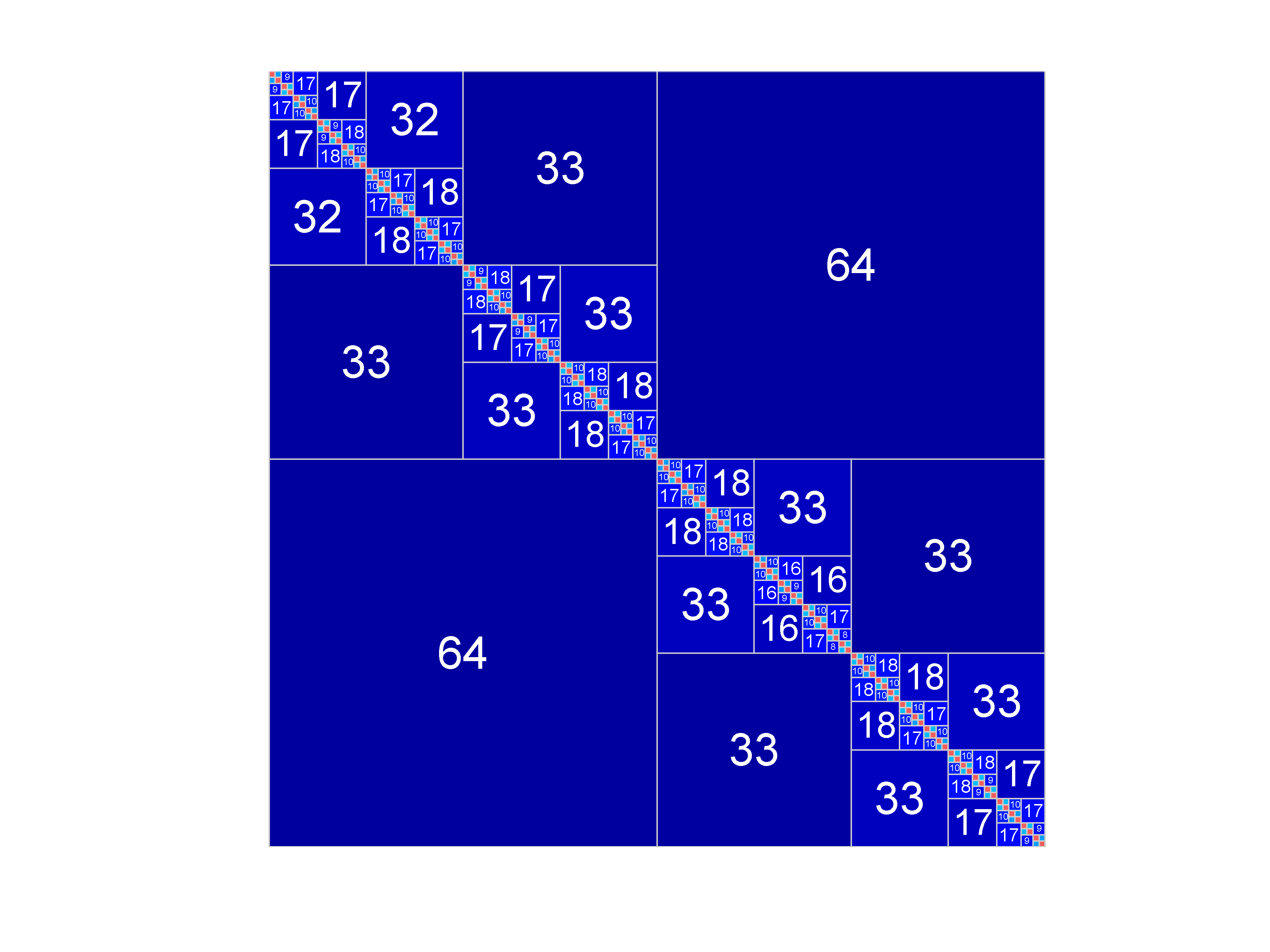}
		\caption{Weak admissibility.}
	\end{subfigure}
	\begin{subfigure}{.5\textwidth}
		\centering
		\includegraphics[width=0.9\linewidth]{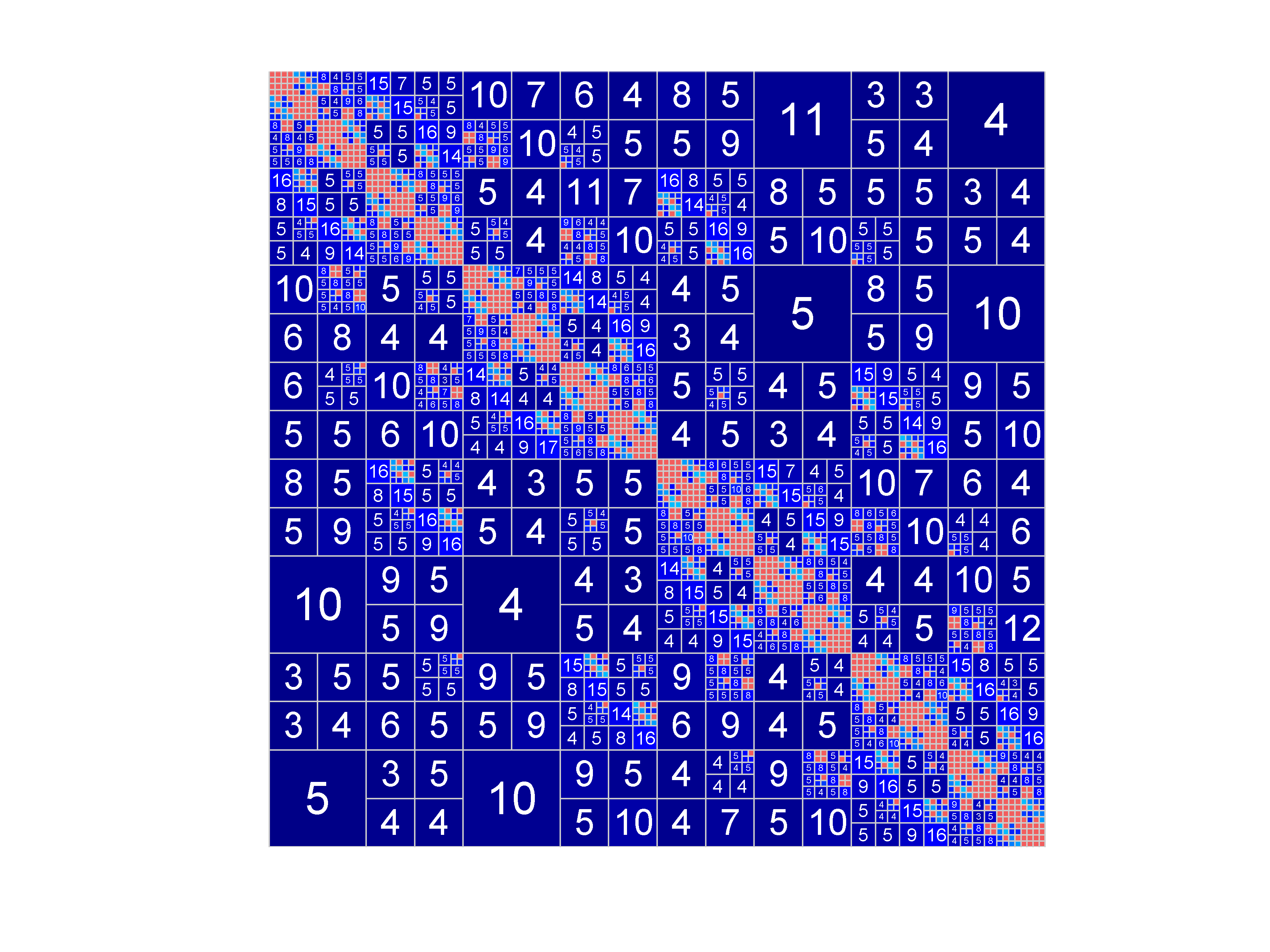}
		\caption{Standard admissibility.}
	\end{subfigure}
	\begin{subfigure}{1.0\textwidth}
		\centering
		\includegraphics[width=0.6\linewidth]{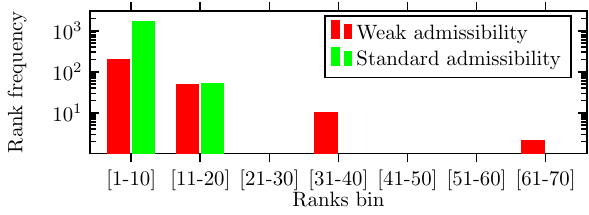}
		\caption{Ranks histogram with respect to admissibility condition.}
	\end{subfigure}
\caption{$\mathcal{H}$-inverse of the 2D Poisson operator, discretized with $N=64\times64$ grid points, using two different admissibility conditions. The number in each block is the numerical rank necessary to achieve a compression accuracy of 1E-4. The color map is determined by the ratio of the numerical rank and the size of the block, deep blue indicates an effective low-rank approximation, while red depict dense blocks.}
\label{fig:admisibilities}
\end{figure}

The motivation for choosing a standard admissibility condition is that, by further refining the off-diagonals blocks, it is possible to achieve a similar accuracy with smaller numerical ranks, that are crucial to ensure economic memory consumption and overall high performance. The impact of the admissibility condition is illustrated in Figure \ref{fig:admisibilities}, which depicts the $\mathcal{H}$-inverse of the variable-coefficient two-dimensional Poisson operator discretized on a $N=64\times64$ grid using a finite difference scheme. In the right panel, the use of a few small dense blocks in the off-diagonal regions allows much smaller ranks to be used in the remaining low-rank blocks, without compromising accuracy.

\nocite{trefethen1997numerical}

Table \ref{table:flopsTable} shows the storage gains by representing the inverse of a 2D Poisson problem with an $\mathcal{H}$-matrix with weak  admissibility versus standard admissibility. Table \ref{table:flopsTable} also shows the difference in terms of number of operations between these two structures. The cost of the $\mathcal{H}$-matrix inversion requires $56 C_{sp}^3 k n(\log n +1)^2 + 184 C_{sp} k^3 n (\log n +1)$ operations, where $k$ represents the average rank of the low-rank blocks, and $C_{sp}$ represents the sparsity of the structure of the hierarchical matrix inverse, see \cite{hackbusch2015hierarchical}. Since the weak admissibility condition requires larger ranks than the standard admissibility condition, at scale, this tends to increase the memory requirements and the number of floating-point operations.

\begin{table}[H]
\centering
\begin{tabular}{|c|c|c|c|c|c|}
\hline
Operation & Format & Storage &  Operations \\ \hline
Inverse & $\mathcal{H}$ (weak admissibility)         & 723 MB  &       8.0E11              \\ \hline
Inverse & $\mathcal{H}$ (standard admissibility)     & 434 MB  &       5.0E11              \\ \hline
Factorization & HSS (weak admissibility)*    & 40 MB  &       3.3E07              \\ \hline
\end{tabular}
\caption{$\mathcal{H}$-inverse of the 2D Poisson operator for $N=256^2$ grid points, using two different admissibility conditions. We document the memory and floating-point operations to build the $\mathcal{H}$-matrix inverse with weak and standard admissibility. The weak admissibility condition tends to require large ranks, which lead to increased memory requirements and more arithmetic operations than the standard admissibility condition. (Equivalent dense storage and arithmetic operations of the inverse operation would have required 2,147 MB and 1.0E13 operations.). *As a matter of comparison we show the equivalent storage requirements of the same problem by using the nested-basis HSS format which for this particular problem has an optimal complexity.}
\label{table:flopsTable}
\end{table}

A low-rank approximation for a given off-diagonal block can be found in a variety of ways. Several strategies, ranging from randomized algorithms to heuristics for pivoting, are available in the literature.  Every block of the $\mathcal{H}$-matrix stored as a low-rank approximation has the form of an outer product $AB^T$. The goal of efficient hierarchical matrix processing is to construct the best possible low-rank factorization as matrix operations are performed. This routine is often referred to as the compression step. For a comprehensive discussion of the construction of $\mathcal{H}$-matrices and its arithmetics, we refer the reader to \cite{hackbusch2015hierarchical}.

\subsection{Cyclic reduction}
\label{cyclic_reduction}

This section reviews the cyclic reduction algorithm in preparation for the following section describing the accelerated cyclic reduction variant that improves its arithmetic and memory complexity growth.

\subsubsection{Model problem}

Consider the seven-point stencil finite difference discretization with Dirichlet boundary conditions of the three-dimensional variable-coefficient Poisson equation on the unit cube.

\begin{equation}
-\nabla \cdot \kappa(x) \nabla u=f(x)
\label{probSta}
\end{equation}

This discretization leads to a block tridiagonal linear system of $N=n^3$ unknowns. This corresponds to a matrix $A$ composed of $3n-2$ blocks of size $n^2 \times n^2$.

\begin{equation}
{\renewcommand{\arraystretch}{1.5}
A = \mbox{tridiagonal}(E_{i},D_{i},F_{i}) =
\begin{bmatrix}
D_0 		& F_0 		&  			&  			&  			\\
E_1 		& D_1 		& F_1 		& 			& 			\\
			& ~\ddots 	& ~\ddots 	& ~\ddots 	& 			\\
			&  			& E_{n-2} 	& D_{n-2} 	& F_{n-2} 	\\
			&  			&  			& E_{n-1} 	& D_{n-1} 	\\
\end{bmatrix}.
}
\label{originalMatrix}
\end{equation}

Block cyclic reduction can be used to solve the system defined in Equation \ref{originalMatrix}. The algorithm consists of two phases: elimination and back substitution. 

\subsubsection{Elimination}

The first step is to rearrange the linear system via matrix permutation $(PAP^T)(Pu)=Pf$. The permutation matrix $P$ corresponds to a red/black (even/odd) ordering of the blocks. For illustration, we choose $n=8$ and consider a $2 \times 2$ partition of the permuted system as shown in Equations \ref{permuted_system} and \ref{partition}. Superscripts indicate step number, where at each step a Schur complementation of a permuted system is performed to reduce the number of unknowns by half.

\begin{equation}
{\arraycolsep=1.4pt\def\arraystretch{1.2}
\begin{array}{cccc}
\left[\begin{array}{cccc|cccc}
D_0^{(0)} 	&     		&     		&     		& F_0^{(0)} &     		&     		&     \\
    		& D_2^{(0)} &     		&     		& E_2^{(0)} & F_2^{(0)} &     		&     \\
    		&     		& D_4^{(0)} &     		&     		& E_4^{(0)} & F_4^{(0)} &     \\
    		&     		&     		& D_6^{(0)} &     		&     		& E_6^{(0)} & F_6^{(0)} \\ \hline 
E_1^{(0)} 	& F_1^{(0)} &     		&     		& D_1^{(0)} &     		&     		&     \\
    		& E_3^{(0)} & F_3^{(0)} &     		&     		& D_3^{(0)} &     		&     \\
    		&     		& E_5^{(0)} & F_5^{(0)} &     		&     		& D_5^{(0)} &     \\
    		&     		&     		& E_7^{(0)} &     		&     		&     		& D_7^{(0)} \\
\end{array}\right]
&
\left[\begin{array}{c}
u_0^{(0)}\\
u_2^{(0)}\\
u_4^{(0)}\\
u_6^{(0)}\\\hline 
u_1^{(0)}\\
u_3^{(0)}\\
u_5^{(0)}\\
u_7^{(0)}\\
\end{array}\right]
&
=
&
\left[\begin{array}{c}
f_0^{(0)}\\
f_2^{(0)}\\
f_4^{(0)}\\
f_6^{(0)}\\\hline 
f_1^{(0)}\\
f_3^{(0)}\\
f_5^{(0)}\\
f_7^{(0)}\\
\end{array}\right]
\end{array}.
}
\label{permuted_system}
\end{equation}

\begin{equation}
{\renewcommand{\arraystretch}{1.8}
\begin{array}{cccc}
\left[\begin{array}{c|c}
A_{11} 	& A_{12} \\ \hline 
A_{21} 	& A_{22} \\\end{array}\right]
&
\left[\begin{array}{c}
u_{even} \\ \hline 
u_{odd} \\\end{array}\right]
&
=
&
\left[\begin{array}{c}
f_{even} \\ \hline 
f_{odd} \\\end{array}\right]
\end{array}.
\label{partition}
}
\end{equation}

The Schur complement computations of the partitioned system are shown in equation \ref{schur_complement}: 
\begin{equation}
(A_{22}-A_{21}A_{11}^{-1}A_{12})u_{odd}=f,\;\;\;\;\;\;\;\; f = f_{odd}-A_{21}A_{11}^{-1}f_{even}.
\label{schur_complement}
\end{equation}

Since the upper-left block $A_{11}$ is block-diagonal, its inverse can be computed as the inverse of each individual block (in this case: $D_0^{(0)}$, $D_2^{(0)}$, $D_4^{(0)}$, and $D_6^{(0)}$), in parallel. All computations for the generation of the Schur complement at step $i+1$, whose size is half of the step $i$ problem, are also done at block-level granularity as show in Equation \ref{elim}, which applies to $j$ odds only. There is a slight abuse of notation in Equation \ref{elim} to handle the case of the last plane that has one neighbor, the computations involving the plane $j+1$ are not performed. We use a polymorphic notation for the matrix addition, matrix subtraction, matrix-matrix multiplication, matrix-vector multiplication, and matrix inversion ($A+_{\mathcal{H}}B, \, A-_{\mathcal{H}}B, \, A\cdot_{\mathcal{H}}B, \, A\cdot_{\mathcal{H}}b, \, A^{-\mathcal{H}}$), depending on whether the matrices are represented in the regular sparse format or the $\mathcal{H}$-matrix format, as we will later refer back when describing the $\mathcal{H}$-matrix accelerated cyclic reduction method.

\begin{equation} 
{\renewcommand{\arraystretch}{1.0}
\begin{array}{l}
E_{j}^{(i+1)}= - E_{j}^{(i)} \cdot_{\mathcal{H}} (D_{j-1}^{(i)})^{-\mathcal{H}} \cdot_{\mathcal{H}} E_{j-1}^{(i)}\\
D_{j}^{(i+1)}=D_{j}^{(i)} -_{\mathcal{H}} \, E_{j}^{(i)} \cdot_{\mathcal{H}} (D_{j-1}^{(i)})^{-\mathcal{H}} \cdot_{\mathcal{H}} F_{j-1}^{(i)} -_{\mathcal{H}} \, F_{j}^{(i)} \cdot_{\mathcal{H}} (D_{j+1}^{(i)})^{-\mathcal{H}} \cdot_{\mathcal{H}} E_{j+1}^{(i)}\\
F_{j}^{(i+1)}= - F_{j}^{(i)} \cdot_{\mathcal{H}} (D_{j+1}^{(i)})^{-\mathcal{H}} \cdot_{\mathcal{H}} F_{j+1}^{(i)}\\\\
f_{j}^{(i+1)}=f_{j}^{(i)} - E_{j}^{(i)} \cdot_{\mathcal{H}} (D_{j-1}^{(i)})^{-\mathcal{H}} \cdot_{\mathcal{H}} f_{j-1}^{(i)} - F_{j}^{(i)} \cdot_{\mathcal{H}} (D_{j+1}^{(i)})^{-\mathcal{H}} \cdot_{\mathcal{H}} f_{j+1}^{(i)}\\
\end{array}
}
\label{elim}
\end{equation}

This process of permuting and Schur complementation is recursive. It finishes when a single block is left, or when the remaining system is small enough to be inverted directly. Recursion is possible because the Schur complement of a tridiagonal matrix is tridiagonal. This property can be seen in the structure of the matrix at the next step shown in Equation \ref{next_step} and illustrating the remaining (originally odd) unknowns after they have been renumbered sequentially.

\begin{equation}
{\renewcommand{\arraystretch}{1.8}
\begin{array}{cccc}
\left[\begin{array}{ccccc}
D_0^{(1)} & F_0^{(1)} &     &     \\
E_1^{(1)} & D_1^{(1)} & F_1^{(1)} &     \\
    & E_2^{(1)} & D_2^{(1)} & F_2^{(1)} \\
    &     & E_3^{(1)} & D_3^{(1)} \\
\end{array}\right]
&
\left[\begin{array}{c}
u_0^{(1)}\\
u_1^{(1)}\\
u_2^{(1)}\\
u_3^{(1)}\\
\end{array}\right]
&
=
&
\left[\begin{array}{c}
f_0^{(1)}\\
f_1^{(1)}\\
f_2^{(1)}\\
f_3^{(1)}\\
\end{array}\right]
\end{array}
}
\label{next_step}
\end{equation}

The algorithm proceeds to apply the red/black permutation followed by a Schur complementation for two more steps to compute the last single block $D_0^{(3)}$. 

\subsubsection{Back-substitution}

Once elimination is completed, the solve stage starts from the last block of unknowns, as shown in equation \ref{first_backSubs}:

\begin{equation}
D_0^{(3)} \cdot_{\mathcal{H}} u_0^{(3)} = f_0^{(3)}.
\label{first_backSubs}
\end{equation}

Once the solution at the last step $u_0^{(3)}$ is computed, it is propagated backward in the hierarchy of the elimination tree.

The formula to compute the solution at step $q$ is given by

\begin{equation}
u^{(i)} = (D^{(i)})^{-\mathcal{H}} \cdot_{\mathcal{H}} ( f^{(i)} - E^{(i)} \cdot_{\mathcal{H}} u^{(i+1)} -  F^{(i)} \cdot_{\mathcal{H}} u^{(i+1)}).
\label{all_backSubs}
\end{equation}

This procedure continues until the solution of the entire linear system is computed.

Back-substitution is much more lightweight than the elimination algorithm regarding computation and communication volume, because it communicates parts of the solution in the form of vectors, and the only matrix operation performed is a matrix-vector multiplication. For large scale problems, this makes the solve phase orders of magnitude faster than the elimination phase. As with other direct solvers, the ability to efficiently solve for a given right-hand side given a factorization motivates the use of ACR for multiple right-hand sides at a minimal cost per new forcing term.

\section{Accelerated Cyclic Reduction}

This section describes how cyclic reduction can be used in combination with hierarchical matrices to result in a variant that improves the computational complexity and memory requirements of the classical cyclic reduction method.

\subsection{Block-wise $\mathcal{H}$-matrix approximation}

ACR approximates each $D_i$, $E_i$ and $F_i$ block of the original block tridiagonal matrix $A$ given in Equation \ref{originalMatrix} with a hierarchical matrix, and then proceeds with the cyclic reduction algorithm, as described in the previous section, by using hierarchical matrix arithmetics instead of the conventional dense linear algebra arithmetic operations.

In generating the structure of the hierarchical matrix representations of the blocks, we exploit the fact the domain is subdivided into $n$ planes each consisting of $n^2$ grid points and block rows of the matrix are identified with the planes of the discretization grid. We consider this geometry and use a two-dimensional planar bisection clustering when constructing each $\mathcal{H}$-matrix. 

Cyclic reduction requires hierarchical matrix addition, subtraction, matrix-matrix multiplication, matrix-vector multiplication and matrix inversion. The relative accuracy of the approximation is specified during the compression of each block and while performing hierarchical matrix arithmetic operations. Committing to a given tolerance ensures that the numerical ranks are adjusted to preserve the specified accuracy during the elimination and solve phases. It is at the block level that the improvements in the complexity estimates take place.

Table \ref{advBlock} summarizes the advantages of a block-wise approximation of matrix blocks with $\mathcal{H}$-matrices in the computation of the inverse of a block, and its storage, as compared to their equivalent dense counterparts.

\begin{table}[H]
\centering
\begin{tabular}{l|l|l}
 					& Inverse 					& Storage\\ \hline
Dense Matrix 		& $\mathcal{O}(N^3)$ 		& $\mathcal{O}(N^2)$ \\
$\mathcal{H}$ Matrix & $O(k~N \log N~(\log N + k^2))$ & $\mathcal{O}(k~N \log N)$ \\
\end{tabular}
\caption{Comparing the complexity estimates of storing and computing the inverse of a $N \times N$ matrix block in dense format, versus approximating the matrix block with a hierarchical matrix with numerical rank $k$.}
\label{advBlock}
\end{table}

\subsection{General algorithm}

To simplify the exposition, we assume the size of the linear system is a power of two; the number of steps required by ACR is thus $q=\log N$. The size of the blocks for 2D problems is $n^2$.

As mentioned in Section \ref{cyclic_reduction}, two procedures define cyclic reduction: elimination and back-substitution.
The high-level algorithm of elimination is shown in listing \ref{acr_elimination}, whereas the high-level algorithm for back-substitution is shown in listing \ref{acr_backsubstitution}. Even though Algorithms \ref{acr_elimination} and \ref{acr_backsubstitution} show permutations and matrix operations at the level of the global system, our implementation operates at a per-block granularity, which means that permutations are part of the implementation's logic and that linear algebraic operations are performed block by block as shown in Equation~\ref{elim}. This is possible since cyclic reduction preserves the block tridiagonal structure during elimination.

\begin{algorithm}[H]
\caption{ACR elimination}
\begin{algorithmic}[1]
\item Block-wise low-rank approximation of $A$: $A^{(0)}$ = $A$
\FOR {$\textit{i}$ = 0 \textbf{to} $\textit{q-1}$}
\item // Generate $A^{(i+1)}$ block tridiagonal of size $n/2^{i+1}$ using block-level operations (Equation \ref{elim})
\item // Requires $O(k~n \log n~(\log n + k^2) / 2^{i+1})$ operations
\STATE $A^{(i+1)}$ = $A_{22}^{(i)} -_{\mathcal{H}} \, A_{21}^{(i)} \cdot_{\mathcal{H}} (A_{11}^{(i)})^{-\mathcal{H}} \cdot_{\mathcal{H}} A_{12}^{(i)}$
\item // Forward substitution, requires $O(k~n \log n/ 2^{i+1})$ operations
\STATE $f^{(i+1)}$ = $f^{(i)}_{2} - A_{21}^{(i)} \cdot_{\mathcal{H}} (A_{11}^{(i)})^{-\mathcal{H}} \cdot_{\mathcal{H}} f^{(i)}_{1}$
\ENDFOR
\end{algorithmic}
\label{acr_elimination}
\end{algorithm}

\begin{algorithm}[H]
\caption{ACR back-substitution}
\begin{algorithmic}[1]
\STATE Solve $A^{(q)} u^{(q)} = f^{(q)}$
\FOR {$\textit{i}$ = $\textit{q-1}$ \textbf{to} 0}
\item // Back-substitution, requires $O(k~n \log n/ 2^{i+1})$ operations
\item // This is performed at block-level (Equation \ref{all_backSubs})
\STATE $u^{(i)} = (A_{11}^{(i)})^{-\mathcal{H}} \cdot_{\mathcal{H}} ( f^{(i)} - A_{12}^{(i)} \cdot_{\mathcal{H}} u^{(i+1)})$
\ENDFOR
\end{algorithmic}
\label{acr_backsubstitution}
\end{algorithm}

\subsection{Sequential complexity estimates}

Every cyclic reduction step requires two matrix-matrix multiplications, one matrix inversion and one matrix addition per block being eliminated. These kernels have arithmetic complexity of $O(k~n \log n~(\log n + k^2))$ operations \cite{hackbusch2015hierarchical}. 
For a problem size of $N = n^3$ with $n = 2^q$, ACR requires $n/2+n/4+n/8+\ldots \approx {n}$ steps to perform elimination. 
The most expensive computation in each step is the computation of an inverse of a block of size $n^2 \times n^2$, which in $\mathcal{H}$-format has a complexity of $O(k~n^2 \log n~(\log n + k^2))$, therefore, ACR results in a $O(k~N \log N~(\log N + k^2))$ overall algorithm, with $\mathcal{O}(k~N \log N)$ memory requirements. Table \ref{arithEstimates} summarizes the complexity estimates of each of the $\mathcal{H}$ matrix operations involved in ACR. Table \ref{complexity} summarizes the complexity estimates of the classical cyclic reduction algorithms without exploitation of equal blocks versus ACR.

\begin{table}[ht!]
\centering
\begin{tabular}{|c|l|}
\hline
Operation & Complexity \\ \hline
$ A +_{\mathcal{H}} B     $ & $O(k^2~n \log n)$  \\ \hline
$ A \cdot_{\mathcal{H}} B $ & $O(k~n \log n~(\log n + k^2))$              \\ \hline
$ A^{-\mathcal{H}}        $ & $O(k~n \log n~(\log n + k^2))$  \\ \hline
$ A \cdot_{\mathcal{H}} b $ & $O(k~n \log n) $              \\ \hline
\end{tabular}
\caption{Summary of the complexity of the $\mathcal{H}$ matrix arithmetic operations.}
\label{arithEstimates}
\end{table}

\begin{table}[H]
\centering
\begin{tabular}{l|l|l}
Method							& Operations 					& Memory \\ \hline
Cyclic Reduction (CR) 				& $\mathcal{O}(N^2)$ 			& $\mathcal{O}(N^{1.5} \log N)$ \\
Accelerated Cyclic Reduction (ACR)    & $O(k~N \log N~(\log N + k^2))$ 	& $\mathcal{O}(k~N \log N)$ \\
\end{tabular}
\caption{Summary of the sequential complexity estimates of the classic cyclic reduction method and the proposed variant, accelerated cyclic reduction; $k$ represents the numerical rank of the approximation.}
\label{complexity}
\end{table}

Because ACR effectively uses hierarchical representations only for a set of regular two-dimensional problems, the resulting constants appearing in the asymptotic complexity estimates tend to be smaller, as a virtue of lower rank requirements, and make it feasible to perform large scale computations. For instance, our limited experiments show that for the 3D Poisson problem (Table \ref{table:ACRPoisConfiguration}) ACR requires substantially lower numerical ranks than the ranks reported in the HSSMF literature \cite{Wang2016,Ghysels15}.

In terms of practical usage, ACR has different concurrency properties than $\mathcal{H}$-LU or multifrontal HSS, enabling different amounts of independent work to be performed. The regularity of the computational patterns of ACR is valuable in terms of the ability to efficiently use current and future hardware architectures.

\section{Parallel accelerated cyclic reduction}
\label{parallelACR}
This section describes how to leverage the concurrency features of the accelerated cyclic reduction method in a distributed-memory parallel environment.

\subsection{Parallel implementation}
The parallel ACR elimination and back-substitution algorithms are listed in Algorithms \ref{ParalllelACRElimination} and \ref{ParalllelACRBacksubstitution}, respectively. 

\begin{algorithm}[H]
\caption{Parallel ACR elimination}
\begin{algorithmic}[1]
\STATE $j$= Processor number

\STATE \textbf{parallel for} at \textit{all} processors $j$, $j \in 0:2^{q}-1$
\STATE \hspace{\algorithmicindent} Block-wise conversion to $\mathcal{H}$-matrix of $A$ = \mbox{tridiagonal}($E_{j}^{(1)},D_{j}^{(1)},F_{j}^{(1)}$)
\STATE \textbf{end parallel for}

\FOR {$\textit{i}$ = 1 \textbf{to} $\textit{q}$}

\STATE \textbf{parallel for} at $j$ \textit{even}, $j \in 0:2^{q-i}-1$
\STATE \hspace{\algorithmicindent} Compute ($D_{j}^{(i)})^{-\mathcal{H}}$
\STATE \hspace{\algorithmicindent} \textbf{Communicate} $E_{j}^{(i)}$, $(D_{j}^{(i)})^{-1}$, $F_{j}^{(i)}$, $f_{j}^{(i)}$ to processors $j-1$ 
\STATE \hspace{\algorithmicindent} \textbf{Communicate} $E_{j}^{(i)}$, $(D_{j}^{(i)})^{-1}$, $F_{j}^{(i)}$, $f_{j}^{(i)}$ to processors $j+1$
\STATE \textbf{end parallel for}

\STATE \textbf{parallel for} at $j$ \textit{odd}, $j \in 0:2^{q-i-1}-1$
\STATE \hspace{\algorithmicindent}  \textbf{Compute} $E_{j}^{(i+1)}$, $D_{j}^{(i+1)}$, $F_{j}^{(i+1)}$, $f_{j}^{(i+1)}$ from Equation \ref{elim}
\STATE \textbf{end parallel for}

\ENDFOR
\end{algorithmic}
\label{ParalllelACRElimination}
\end{algorithm}

\begin{algorithm}[H]
\caption{Parallel ACR back-substitution}
\begin{algorithmic}[1]
\STATE $n=2^q$
\STATE $j$= Processor number
\FOR {$\textit{i}$ = $\textit{q}$ \textbf{to} 1}
\STATE \textbf{parallel for} at $j$, $j \in 0:2^{q-i}-1$
\STATE \hspace{\algorithmicindent} \textbf{Compute} $u_{j}^{(i)}$ from Equation \ref{all_backSubs}
\STATE \hspace{\algorithmicindent} \textbf{Communicate} $u_{j}^{(i)}$ to processors $j-1$ 
\STATE \hspace{\algorithmicindent} \textbf{Communicate} $u_{j}^{(i)}$ to processors $j-1$ 
\STATE \textbf{end parallel for}
\ENDFOR
\end{algorithmic}
\label{ParalllelACRBacksubstitution}
\end{algorithm}

A number of concurrency features of the algorithms are evident. Each block row, identified by a plane in the discretization, is assigned to an MPI rank. This decomposition allows the initial conversion of each block into an $\mathcal{H}$-matrix in an embarrassingly parallel manner. The $q = \log n$ levels of Schur complement computation exploit concurrent execution in two ways:

\begin{itemize}
  \item The inverse of the block $A_{11}$ of Equation~\ref{partition} can be computed concurrently in a block-wise fashion since $A_{11}$ is block diagonal. This computation is embarrassingly parallel.
  \item Computing the Schur complement requires two matrix-matrix multiplications and one matrix addition. Since the linear system partition is formed out of matrix blocks, the computation of these block matrix-matrix multiplications and block matrix-addition can also be computed concurrently.
\end{itemize}

\begin{figure}[H]
\centering
\includegraphics[width=0.6\textwidth]{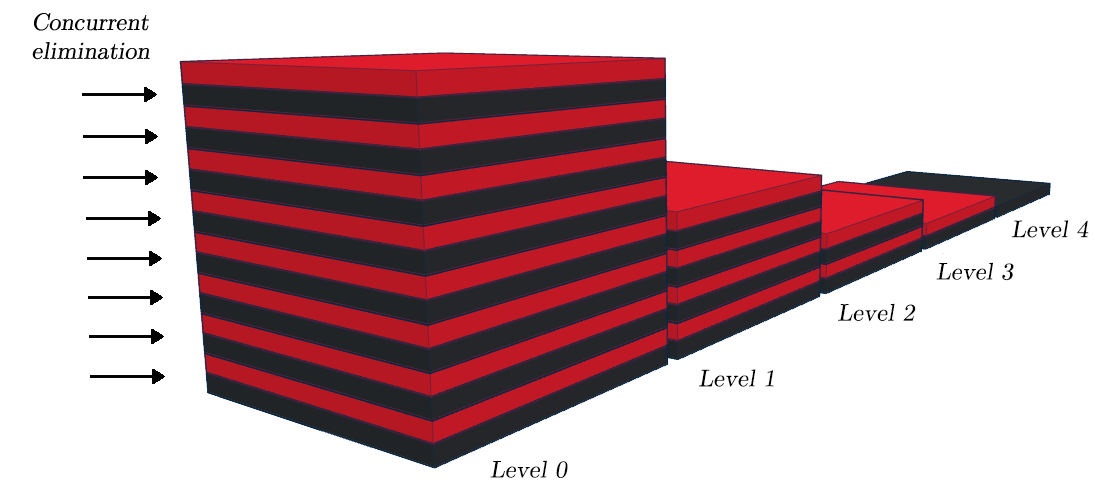}
\caption{Concurrency in ACR elimination for the 16-planes case. Level 0 can eliminate eight planes concurrently thus reducing the problem size to the next level by two; this process continues until one plane is left.}
\label{fig:concurrency}
\end{figure}

Figure~\ref{fig:concurrency} depicts the concurrency through the various levels in ACR elimination. We note here that the ACR decomposition strategy bears a similarity to the slice decomposition \cite{gugist2001}, and also relate to the sweeping preconditioner strategy \cite{engquist2011sweeping}, with the key distinction being that rather than sweeping through the domain, ACR eliminates several planes at once, concurrently. 

\subsection{Inter-node communication}

\begin{figure}[H]
\centering
\includegraphics[width=0.8\textwidth]{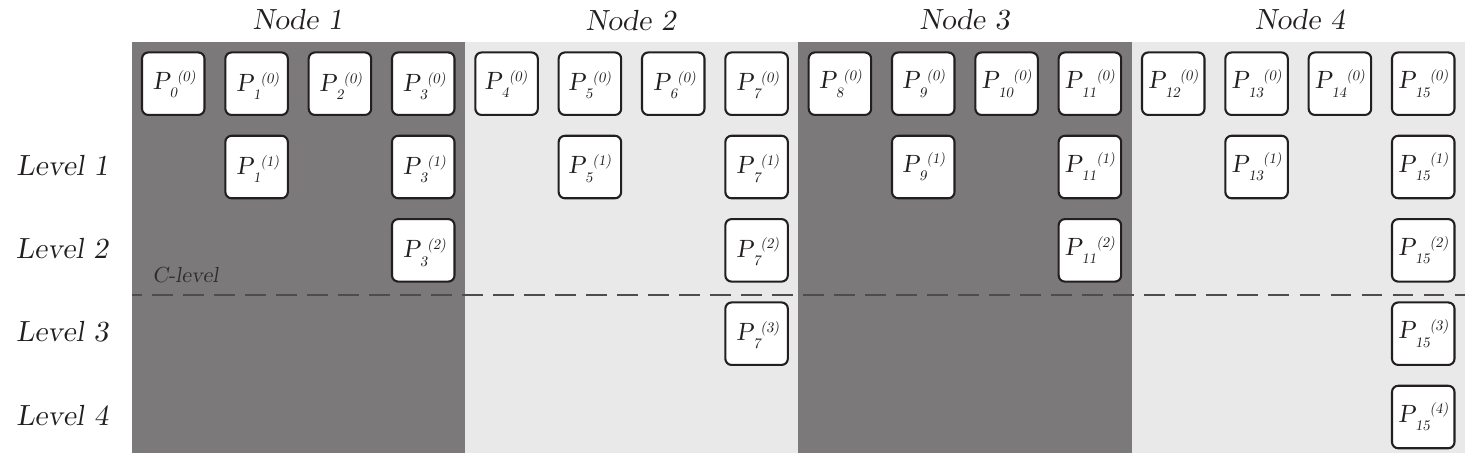}
\caption{Distribution of multiple planes per physical compute node for an example with n=16 and p = 4.}
\label{planes_per_node}
\end{figure}

In the current implementation, each plane is assigned to an MPI rank, and multiple planes are assigned to compute nodes. Let $p$ be the number of physical compute nodes each storing $n/p$ planes at the beginning of the factorization. After $r$ steps of ACR, each compute node holds $n/(2^r p)$ planes. At level $r = \log (n/p)$, a coarse level called the C-level, every node holds a single plane only. The remaining  $\log p$ steps of ACR beyond the C-level leave some compute nodes idle as illustrated in Figure~\ref{planes_per_node}.

\begin{figure}[H]
\includegraphics[width=0.6\textwidth]{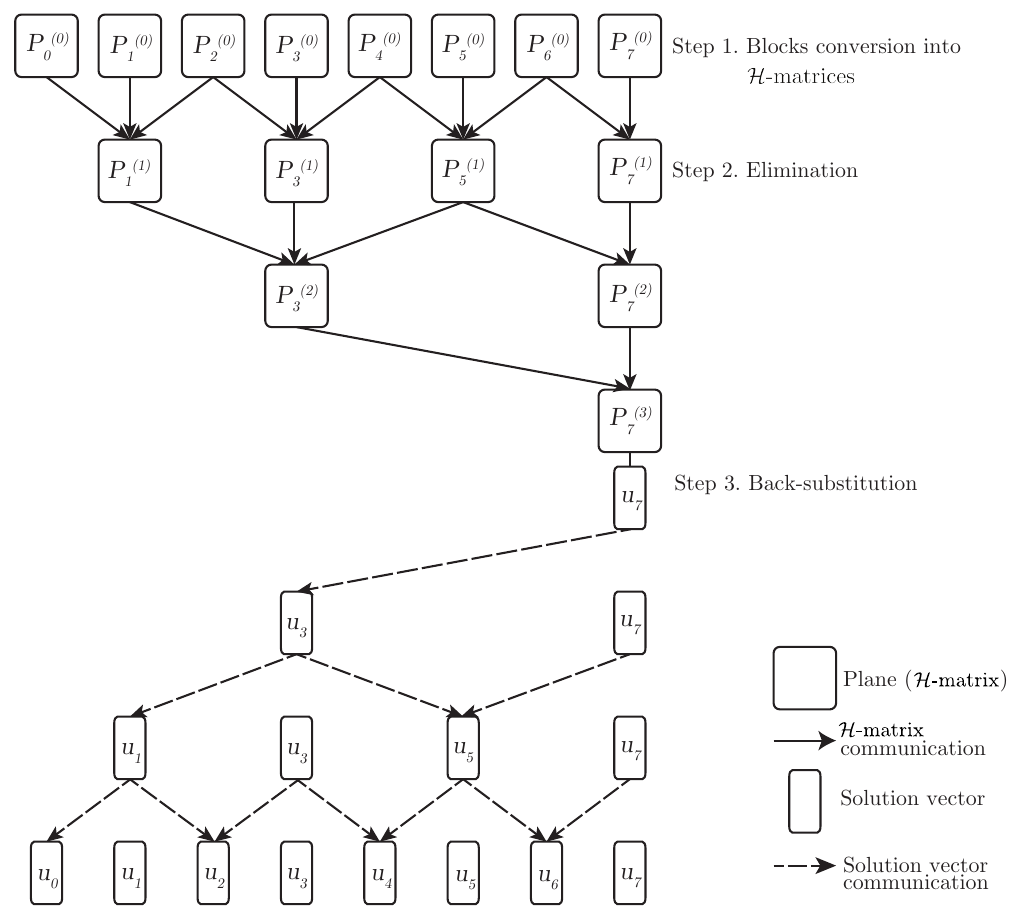}
\centering
\caption{Communication pattern for the 8-planes case. $P_j$ depicts the plane being eliminated, and $u_j$ its corresponding solution.}
\label{CommunicationPattern}
\end{figure}

Distributed-memory communication occurs just at inter-node boundaries thanks to sorting at every step of the factorization, as the computation of the Schur complement for plane $P_{j}$ just requires planes $P_{j-1}$ and $P_{j+1}$, see Figure \ref{planes_per_node}. Thus up to the C-level there are $O(p)$ communication messages per step, each transmitting planes of size $O(k~n^2 \log n)$. Beyond the C-level, there are $O(p/2+\cdots+1) \approx O(p)$ communications messages, adding up to a total communication volume of $O(k~p~n^2 \log n~ ( \log \frac{n}{p}+1 ) )$ for ACR. The communication pattern with its bottom-up binary tree structure is depicted in Fig. \ref{CommunicationPattern}.

\subsection{Parallel time complexity}

The regularity of the ACR algorithm makes it straightforward to estimate the parallel time of the factorization and assess its scalability characteristics. Consider the longest computing node which executes $\log n$ ACR steps. In the $\log (n/p)$ steps preceding the C-level, this node processes $n/(2p) + n/(4p) + \cdots + 1$ block rows in sequence. Beyond the C-level, it processes a single block row in every one of the sequential $\log p$ steps. This results in an asymptotic parallel time complexity for ACR of $O\left( k n^2\log n (\log n + k^2) ( n/p + \log p ) \right)$. The sequential computational time gets reduced by the number of parallel compute nodes $p$, but at the expense of an additional $\log p$ factor that inhibits perfect strong scaling. Fortunately, the amount of work above the C-level that introduces this $\log p$ factor left is small and grows only as $n^2 = N^{2/3}$. 

\begin{figure}[H]
\centering
\includegraphics[width=0.6\textwidth]{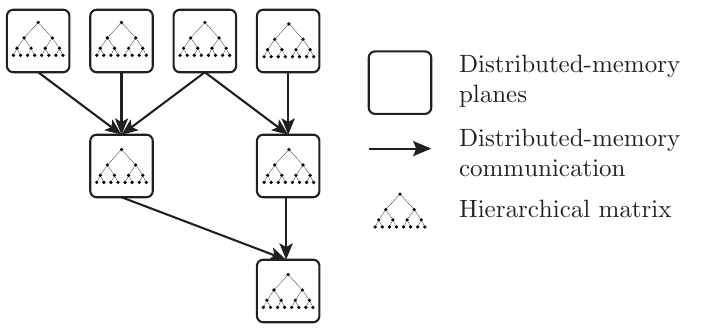}
\caption{Parallel ACR elimination tree depicting two-levels of concurrency using distributed-memory parallelism to distribute concurrent work across compute nodes, and shared-memory parallelism to perform $\mathcal{H}$-matrix operations within the nodes.}
\label{two_levels_parallelism}
\end{figure}

Finally, we note that beyond the parallelism across distributed computing nodes, there is additional concurrency available at the node level. This additional level of parallelism is possible, not only because elimination and back-substitution for multiple block rows can proceed concurrently, but also because parallel variants of the hierarchical matrix arithmetics can be used in performing operations on individual blocks. The two levels of intra-node parallelism are shown schematically in Figure \ref{two_levels_parallelism}. In practice, programming models based on tasks and directed acyclic graphs have proven to be effective to parallelize hierarchical matrix arithmetics \cite{Kriemann2014,Ghysels15}, but the optimal allocation of the multiple cores of a compute node to either block row processing or to individual operations on single blocks requires tuning. We do not describe this aspect of the parallel implementation further here.

\section{Numerical results}
\label{numerical_results}

This section documents the parallel performance and scalability of ACR in a distributed-memory environment. 
The source code is written in C and compiled with the Intel C compiler v15. External libraries utilized in the reference implementation include HLIBpro v2.2 with Intel TBB \cite{kriem05,grasedyck2008performance}, and the sequential version of the Intel Math Kernel Library \cite{kalinkin2015schur}.
Experiments are conducted on the Cray XC40 Shaheen supercomputer at the King Abdullah University of Science $\&$ Technology. Each node has 128GB of RAM and two Intel Haswell processors, each with 16 cores clocked at 2.3Ghz.

To provide a baseline of performance we consider the solution of the same linear systems with STRUMPACK \cite{Ghysels15} v1.0.3, the open-source implementation of the HSS-structured multifrontal solver (HSSMF) developed at the Lawrence Berkeley National Laboratory. The HSSMF method can solve a broader class of linear systems compared to ACR, but the comparison is still of interest, as STRUMPACK is among the few available implementations of distributed-memory fast direct solvers that exploit hierarchically low-rank approximations.

The tuning parameters of ACR include the choice of the leaf node size $n_{min}$ for the $\mathcal{H}$ matrices, the threshold parameter $\eta$ used to decide which blocks will be approximated with a low-rank factorization, or as a dense, full-rank, block, and the accuracy of the approximation $\mathcal{H}_{\epsilon}$ for the construction and algebraic operations of the $\mathcal{H}$ matrices. The tuning parameters for STRUMPACK include how many matrices from the nested-dissection elimination tree will be approximated as HSS, which is controlled by specifying the threshold at which frontal matrices will represented as HSS matrices, the compression accuracy for the HSS matrices, and the minimum leaf size of the HSS frontal matrices. We recall here that the HSS matrix format uses the so-called weak admissibility condition, whereas ACR uses a standard admissibility condition, which does not limit the use of dense blocks exclusively at the matrix diagonal. Additionally, we also consider the algebraic multigrid (AMG) implementation of hypre \cite{Briggs00, hypre_web_page}. Comparison experiments are set to deliver a solution with a relative error tolerance as ${\left|\left|Ax-b\right|\right|_2}/{\left|\left| b \right|\right|_2}\approx~10^{-2}$. For further comparisons, we also consider the multifrontal (MF) implementation of STRUMPACK, and our cyclic reduction (CR) implementation with dense matrix blocks.

\begin{figure}[p]
	\centering
	\begin{subfigure}{.4\textwidth}
		\centering
		\includegraphics[width=0.9\linewidth]{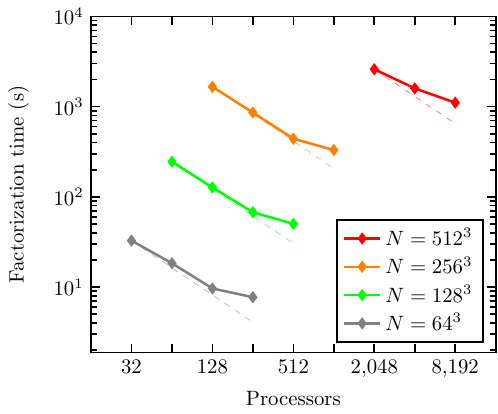}
		\caption{Strong scaling of factorization.}
		\label{fig:pois1}
	\end{subfigure}
	\begin{subfigure}{.4\textwidth}
		\centering
		\includegraphics[width=0.9\linewidth]{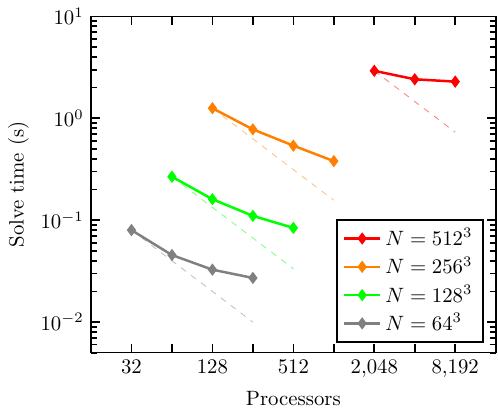}
		\caption{Strong scaling of solve.}
		\label{fig:pois2}
	\end{subfigure}
	\begin{subfigure}{.4\textwidth}
		\centering
		\includegraphics[width=0.9\linewidth]{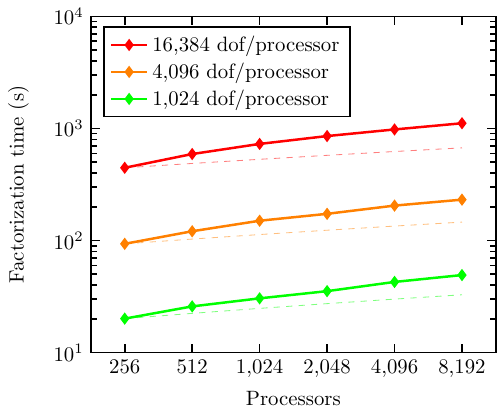}
		\caption{Weak scaling of factorization.}
		\label{fig:pois3}
	\end{subfigure}
	\begin{subfigure}{.4\textwidth}
		\centering
		\includegraphics[width=0.9\linewidth]{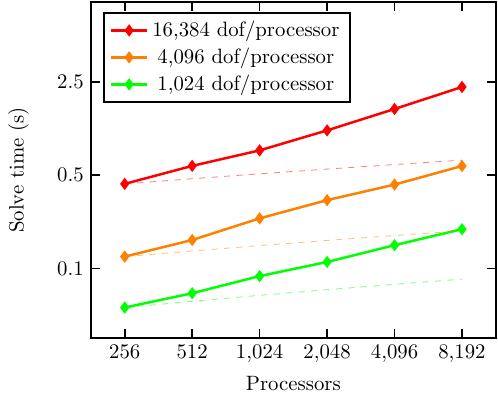}
		\caption{Weak scaling of solve.}
		\label{fig:pois4}
	\end{subfigure}
	\begin{subfigure}{.4\textwidth}
		\centering
		\includegraphics[width=0.9\linewidth]{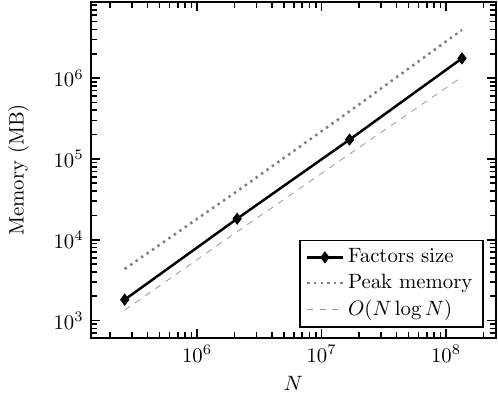}
		\caption{Memory consumption.}
		\label{fig:pois5}
	\end{subfigure}
	\begin{subfigure}{.4\textwidth}
		\centering
		\includegraphics[height=1.6in]{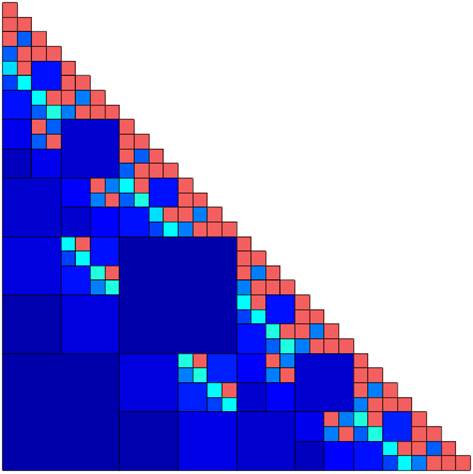}
		\caption{Choice of $\mathcal{H}$-matrix structure to represent planes. Blue indicates low-rank blocks, whereas red indicates dense blocks.}
		\label{fig:pois6}
	\end{subfigure}
\caption{Parallel scalability and memory consumption of ACR for the solution of the Poisson equation. Figure \ref{fig:pois6} depicts the structure of the $\mathcal{H}$-matrices utilized for all the ACR blocks (in this case, extracted from the $N=32^3$ problem). The deep blue color indicates an effective compression, while red blocks indicate no compression.}
\label{fig:scalingPois}
\end{figure}

\subsection{Poisson equation}
\label{section:poisson}

We consider a sequence of Poisson problems of up to $512^3 \approx 134$M unknowns, which is considered very large for this type of ``direct'' (as opposed to iterative) methods. We feature the Poisson equation with homogeneous Dirichlet boundary conditions in the unit cube, i.e.

\begin{equation}
\begin{aligned}
-\nabla^{2} u  = 1, \;\;\; &\mathbf{x} \in \Omega = [0,1]^3,  \;\;\; u(\mathbf{x}) = 0, \; x \in \Gamma, \\
\end{aligned}
\label{poissonEquation}
\end{equation}
discretized with the 7-point finite-difference star stencil, which leads to a symmetric positive definite linear system.

Although this problem can be solved with ACR, other methods such as multigrid or FFTs are ordinarily used instead; we consider it to report on a standard and well-known problem, to facilitate the exposition of ACR. Furthermore, the discretization of the Poisson equation has all positive eigenvalues with rapid decay in off diagonal, making it also an ideal case for hierarchically low-rank approximations analysis.

Figures \ref{fig:pois1} and \ref{fig:pois2} show the total time in seconds for the factorization and solve phases of ACR in a strong scaling setting; dashed lines indicate ideal scaling. Ideal scaling of the factorization stage deteriorates at large processor counts as factors such as communication volume and hardware latency begin to play a significant role; the same factors tend to dominate even more during the solve phase, being the latter a sequence of fast $\mathcal{H}$-matrix-vector multiplications with limited availability of communication/computation overlap. 

Figures \ref{fig:pois3} and \ref{fig:pois4} depict the results of a weak scaling test for ACR with different numbers of degrees of freedom per processor, along with the ideal weak scaling reference lines depicted as dashed curves. The timings deviate from the ideal scaling due to the inherently load imbalance of the recursive bisection strategy of cyclic reduction as some processors become idle towards the end of the reduction. Communication latency further impacts the solve stage at large core counts due to the lower arithmetic intensity of this stage.

Figure \ref{fig:pois5} depicts the memory requirements to store the ACR factorization, together with the expected asymptotic memory usage as $O(N \log N)$.
We stress that the maximum rank of the factored matrices varies from 5 to 10 within all the combinations of problem sizes/number of processors considered in the strong and weak scaling tests (data not shown). Figure \ref{fig:pois6} depicts the structure of the $\mathcal{H}$-matrices used to represent each plane, with the choice of standard admissibility condition. Dark blue blocks denote a low ratio between the numerical rank of the approximation and the full rank of the block, whereas red block indicates non-admissible blocks stored in dense format. For visualization purposes, the figure was taken from the $N=32^3$ problem, and represents the last diagonal block during the elimination phase of ACR. The prevalence of dark blue blocks indicate a good relative compression of each block, since the ratio of numerical rank of the approximation and the actual block size is very small. Most of the red blocks are clustered near the diagonal, where the smallest blocks reside.

\begin{figure}[]
	\centering
	\begin{subfigure}{.4\textwidth}
		\centering
		\includegraphics[width=0.8\linewidth]{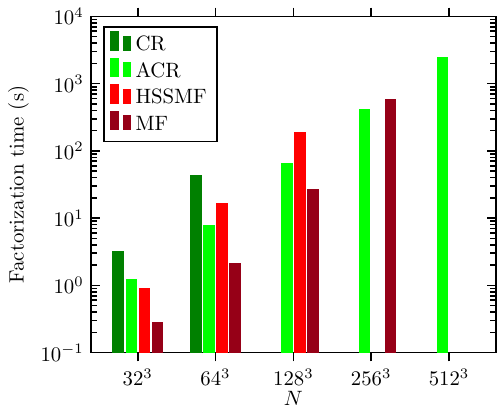}
		\caption{Factorization time.}
		\label{fig:pois11}
	\end{subfigure}
	\begin{subfigure}{.4\textwidth}
		\centering
		\includegraphics[width=0.8\linewidth]{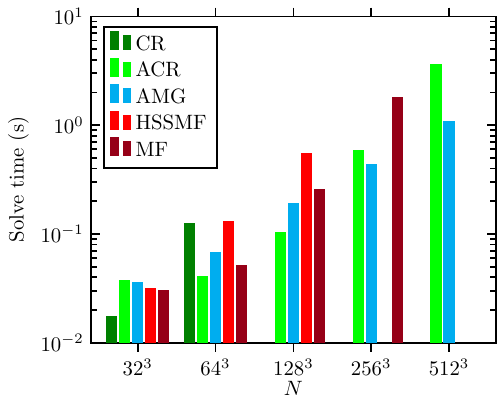}
		\caption{Solve time.}
		\label{fig:pois12}
	\end{subfigure}
	\begin{subfigure}{.4\textwidth}
		\centering
		\includegraphics[width=0.8\linewidth]{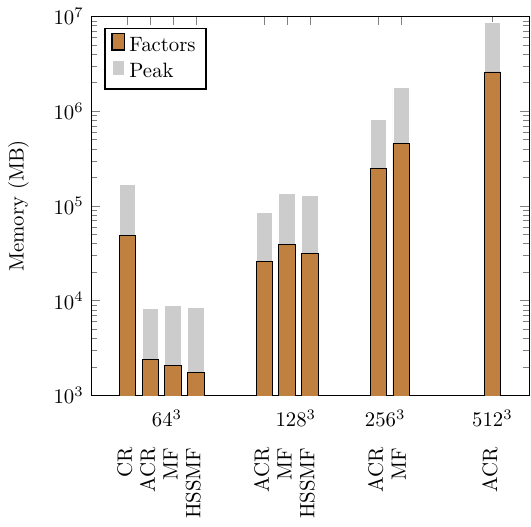}
		\caption{Memory consumption.}
		\label{fig:pois13}
	\end{subfigure}
	\begin{subfigure}{.4\textwidth}
		\centering
		\includegraphics[width=0.8\linewidth]{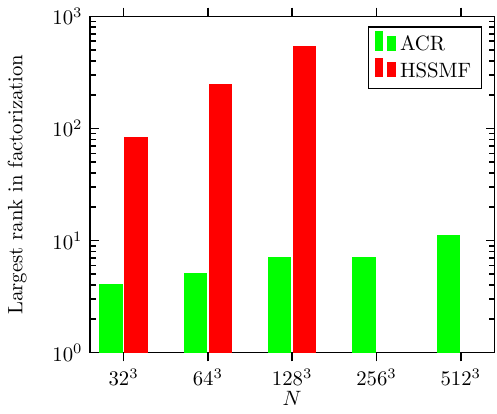}
		\caption{Largest rank in the factorization.}
		\label{fig:pois14}
	\end{subfigure}
\caption{Performance of the factorization and solve phases of ACR for the Poisson problem.}
\label{fig:PoissComparisons}
\end{figure}

Figure \ref{fig:PoissComparisons} compares all solvers under consideration for a set of Poisson problems that range from $N=32^3$ to $N=512^3$ unknowns, with processor counts increased from 256 to 4,096 respectively.
We document the execution parameters, obtained relative residual, and ranks of the ACR and HSSMF factorization in Tables \ref{table:ACRPoisConfiguration} and \ref{table:STRPoisConfiguration}. 
We report factorization times in Figures \ref{fig:pois11} showing that ACR can competitively tackle these problems. Similarly, the solve timings in Figure \ref{fig:pois12} show that ACR is able to solve for a given right-hand size in comparable times as the other methods under consideration. Figure \ref{fig:pois13} documents the size of the factors required by the factorizations, and it shows that the cyclic reduction method (CR) cannot solve problems as small as $N=128^3$ due to memory limitations.
Additionally, we report the peak memory utilization of each solver using the library PAPI v5.5 \cite{browne2000PAPI}, which shows the largest memory usage that each solver required to produce the factorization.
Also, the experiments confirm that the HSSMF method requires less memory to store its factors than the multifrontal method (MF). However, as Figure \ref{fig:pois14} shows, the HSSMF method requires higher ranks than ACR, which translated into a larger size of the factors and prohibited the execution of HSSMF for problems of $N=256^3$ and above. 
The experiments show that ACR requires only $O(1)$ ranks, as opposed to the $O(n)$ rank requirements of the HSSMF factorization.

\begin{table}[]
\centering
\begin{tabular}{|c|c|c|c|c|c|c|}
\hline
\textbf{$N$} & $\mathcal{H}_{\epsilon}$ & \textbf{$\eta$} & \textbf{\begin{tabular}[c]{@{}c@{}}Leaf\\ size\end{tabular}} & \textbf{\begin{tabular}[c]{@{}c@{}}Relative\\ residual\end{tabular}} & \multicolumn{1}{l|}{\textbf{Average rank}} & \multicolumn{1}{l|}{\textbf{Largest rank}} \\ \hline
$32^3$       & 8E-03   & 2               & 32                                                           & 1.39E-02                                                             & 3                                          & 4                                          \\ \hline
$64^3$       & 1E-03   & 2               & 32                                                           & 3.20E-02                                                             & 3                                          & 5                                          \\ \hline
$128^3$      & 1E-03   & 2               & 32                                                           & 2.22E-02                                                             & 4                                          & 7                                          \\ \hline
$256^3$      & 1E-03   & 2               & 32                                                           & 8.75E-02                                                             & 4                                          & 7                                          \\ \hline
$512^3$      & 1E-04   & 2               & 32                                                           & 3.26E-02                                                             & 5                                          & 11                                         \\ \hline
\end{tabular}
\caption{Execution parameters, obtained relative residual, and ranks of the ACR factorization for the Poisson experiments.}
\label{table:ACRPoisConfiguration}

\bigskip

\centering
\begin{tabular}{|c|c|c|c|c|c|c|}
\hline
\textbf{$N$} & \textbf{\begin{tabular}[c]{@{}c@{}}Compression\\ tolerance\end{tabular}} & \multicolumn{1}{l|}{\textbf{\begin{tabular}[c]{@{}l@{}}Relative\\ residual\end{tabular}}} & \textbf{\begin{tabular}[c]{@{}c@{}}Leaf\\ size\end{tabular}} & \textbf{\begin{tabular}[c]{@{}c@{}}Minimum\\ front size\end{tabular}} & \multicolumn{1}{l|}{\textbf{\begin{tabular}[c]{@{}l@{}}Largest rank\end{tabular}}} \\ \hline
$32^3$       & 1E-02                                                                               & 4.41E-02                                                                                  & 128                                                          & 256                                                                   & 82                                                                                                    \\ \hline
$64^3$       & 1E-03                                                                               & 2.65E-02                                                                                  & 128                                                          & 1,024                                                                  & 243                                                                                                   \\ \hline
$128^3$      & 1E-03                                                                               & 8.40E-02                                                                                  & 128                                                          & 4,096                                                                  & 532                                                                                                   \\ \hline
\end{tabular}
\caption{Execution parameters, obtained relative residual, and ranks of the HSSMF factorization for the Poisson experiments.}
\label{table:STRPoisConfiguration}

\end{table}

As expected for this particular problem, multigrid is the method of choice concerning performance and memory footprint for a single right-hand-side. However, for multiple right-hand-sides, the ability to reuse the factorization could give the advantage to solvers based on factorization. The factorization times for ACR and HSSMF are comparable, with the setup stage of HSSMF being faster for smaller problems; the smaller ranks required by ACR instead lead to a faster factorization step with large problem sizes and faster time to solution.

\begin{figure}[]
	\centering
	\begin{subfigure}{.4\textwidth}
		\centering
		\includegraphics[width=0.8\linewidth]{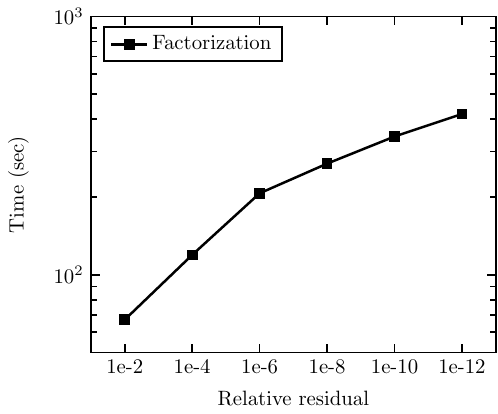}
		\caption{ACR factorization.}
		\label{fig:directFact}
	\end{subfigure}
	\begin{subfigure}{.4\textwidth}
		\centering
		\includegraphics[width=0.8\linewidth]{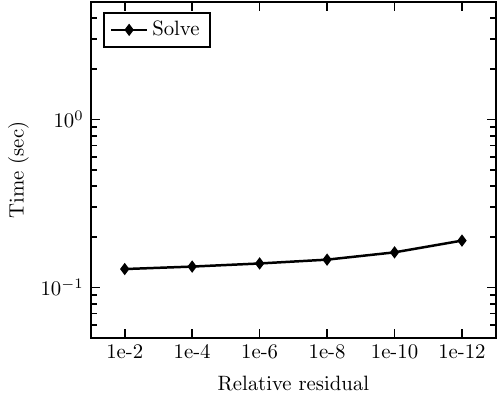}
		\caption{ACR solve.}
		\label{fig:directSolve}
	\end{subfigure}
	\begin{subfigure}{0.4\textwidth}
		\centering
		\includegraphics[width=0.8\linewidth]{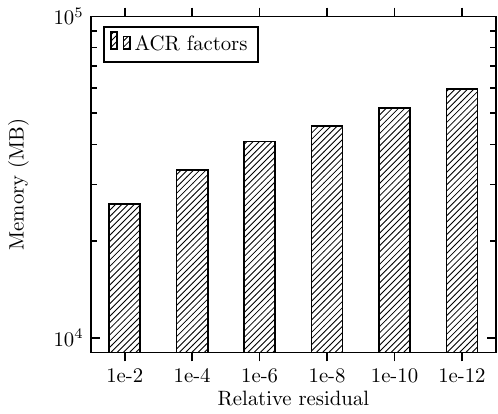}
		\caption{ACR size of the factors.}
		\label{fig:directMem}
	\end{subfigure}
		\begin{subfigure}{.4\textwidth}
		\centering
		\includegraphics[width=0.8\linewidth]{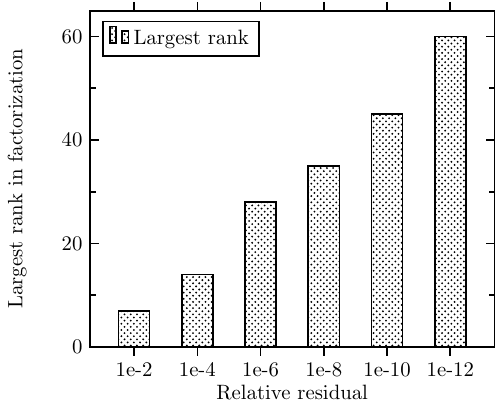}
		\caption{Ranks requirements of the factorization.}
		\label{fig:directRank}
	\end{subfigure}
\caption{Controllable accuracy solution of ACR for a $N=128^3$ Poisson problem.}
\label{fig:controlAccuracy} 
\end{figure}

While ACR and HSSMF solvers can deliver a more accurate solution as direct solvers (i.e. without iterative procedures), this comes at the expense of more time and memory; it is common practice that this factorization is then used as a preconditioner or passed to an iterative refinement procedure. Numerical experiments confirm that ACR could be used as a direct solver if we tune its parameters with a higher accuracy for its $\mathcal{H}$-matrix representations and operations, as depicted in Figure \ref{fig:controlAccuracy}, at the expense of modest rank increases, albeit with higher memory requirements and time to solution. However, as Table \ref{table:precond} shows, a low-accuracy factorization in combination with an iterative procedure is best to minimize the total time-to-solution.

\begin{table}[]
\centering
\begin{tabular}{|l|c|c|c|c|c|c|c|}
\hline
\multicolumn{1}{|c|}{$\mathcal{H}_{\epsilon}$} & \textbf{Factors (MB)} & \textbf{Largest rank} & \textbf{Factorization} & \textbf{Apply} & \textbf{Total time} & \textbf{Iterations} \\ \hline
6E-01                              & 17,280                & 31                 & 18.55                  & 0.050          & 20.72               & 43                  \\ \hline
3E-01                              & 19,385                & 31                 & 21.33                  & 0.053          & 23.14               & 34                  \\ \hline
1E-01                              & 22,328                & 31                 & 26.56                  & 0.058          & 28.01               & 25                  \\ \hline
1E-02                              & 26,687                & 37                 & 51.24                  & 0.064          & 51.94               & 11                  \\ \hline
1E-03                              & 32,212                & 53                 & 89.32                  & 0.104          & 89.73               & 4                   \\ \hline
1E-04                              & 39,181                & 71                 & 149.06                 & 0.127          & 149.44              & 3                   \\ \hline
\end{tabular}
\caption{Iterative solution of a $N=128^3$ Poisson problem with the conjugate gradients method and ACR preconditioner. Relative residual of the solution is 1E-6 in all cases.}
\label{table:precond}
\end{table}

\subsection{Convection-diffusion equation}

We next consider a standard convection-diffusion problem
\begin{equation}
	\begin{aligned}
		&-\nabla^{2} u  + \alpha b(\mathbf{x}) \cdot \nabla u = f(\mathbf{x}), \;\;\; \mathbf{x} \in \Omega = [0,1]^3, \\
		&b(\mathbf{x}) = \begin{bmatrix}
			\sin(a\;2\pi x) \sin(a\;2\pi (1/8 + y)) +  \sin(a\;2\pi (1/8 + z)) \sin(a\;2\pi x) \\
			\cos(a\;2\pi x) \cos(a\;2\pi (1/8 + y)) +  \cos(a\;2\pi (1/8 + y)) \cos(a\;2\pi z) \\
			\cos(a\;2\pi x) \cos(a\;2\pi (1/8 + z)) +  \sin(a\;2\pi (1/8 + y)) \sin(a\;2\pi z)
		\end{bmatrix},\\
	\end{aligned}
\label{eq:codi}
\end{equation}
discretized with a 7-point upwind finite difference scheme, that leads to a non-symmetric linear system which is challenging for classical iterative solvers, especially when the convection term dominates the equation. The $b(\mathbf{x})$ term we consider is a three-dimensional generalization of the two-dimensional vortex flow proposed by Wessel et. al. \cite{ames2014numerical}. We adjust the forcing term and boundary conditions to meet the exact solution
$$u(\mathbf{x}) = \sin(\pi x)+\sin(\pi y)+\sin(\pi z)+\sin(3\pi x)+\sin(3\pi y)+\sin(3\pi z),$$
as proposed by Gupta and Zhang \cite{gupta2000high}, as it is an archetypal challenging problem for multigrid methods.

\begin{figure}[]
	\centering
	\begin{subfigure}{0.4\textwidth}
		\centering
		\includegraphics[width=0.8\linewidth]{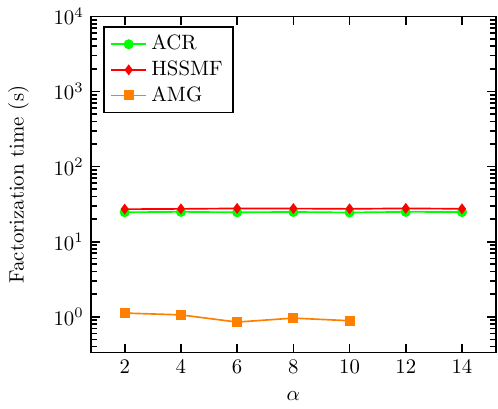}
		\caption{Factorization time.}
		\label{fig:codi1}
	\end{subfigure}
	\begin{subfigure}{0.4\textwidth}
		\centering
		\includegraphics[width=0.8\linewidth]{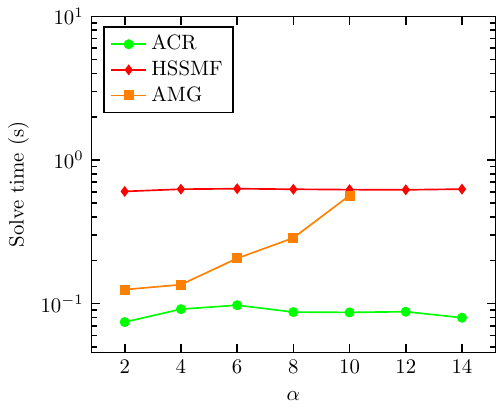}
		\caption{Solve time.}
		\label{fig:codi2}
	\end{subfigure}
	\begin{subfigure}{0.4\textwidth}
		\centering
		\includegraphics[width=0.8\linewidth]{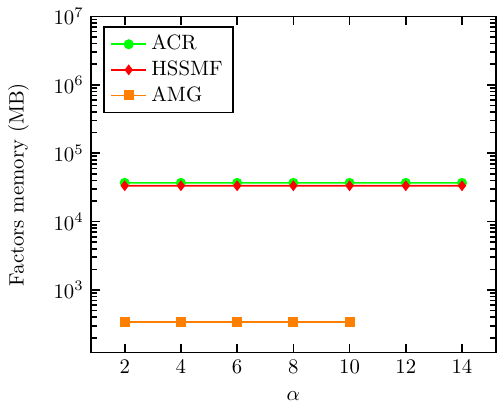}
		\caption{Factors size.}
		\label{fig:codi3}
	\end{subfigure}
	\begin{subfigure}{0.4\textwidth}
		\centering
		\includegraphics[width=0.8\linewidth]{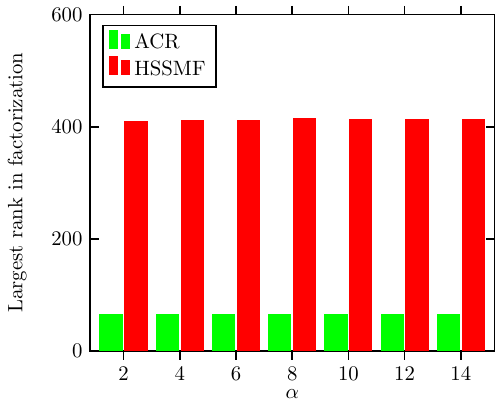}
		\caption{Largest rank in the factorization.}
		\label{fig:codi4}
	\end{subfigure}
\caption{Robustness of ACR for convection-diffusion problem. In convection dominated problems (large $\alpha$), AMG fails to converge while direct solvers maintain a steady performance.}
\label{fig_codi}
\end{figure}

To demonstrate the robustness of ACR and HSSMF for this problem, we fix the number of degrees of freedom at $N=128^3$ and we increase the dominance of the convective term; results are reported in Figure \ref{fig_codi}. Consistently with the Poisson problem, multigrid methods remains the method of choice for diffusion dominated problems in terms of time to solution; however, when $\alpha$ is increased, the performance of AMG deteriorates. On the other hand, both ACR and HSSMF prove to be able to solve convection-dominated problems, with ACR being consistently faster than HSSMF particularly in the back-substitution phase. The size of the factors generated by ACR and HSSMF are comparable, with ACR using significantly smaller ranks.

\subsection{Helmholtz equation}

\begin{figure}[]
	\centering
	\begin{subfigure}{0.4\textwidth}
		\centering
		\includegraphics[width=\linewidth]{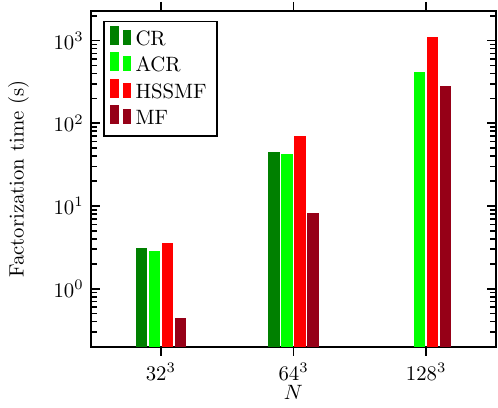}
		\caption{Factorization time.}
		\label{fig:hel1}
	\end{subfigure}
	\begin{subfigure}{0.4\textwidth}
		\centering
		\includegraphics[width=\linewidth]{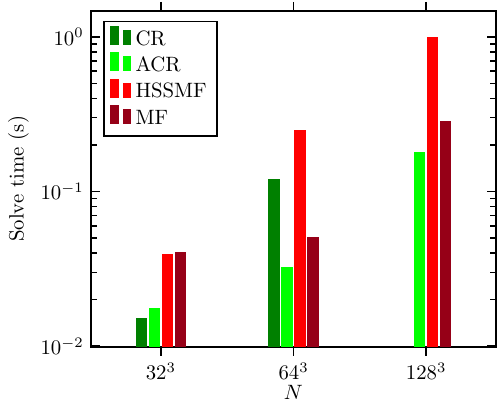}
		\caption{Solve time.}
		\label{fig:hel2}
	\end{subfigure}
	\begin{subfigure}{0.4\textwidth}
		\centering
		\includegraphics[width=\linewidth]{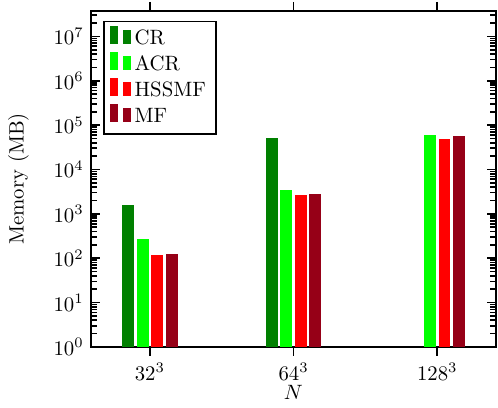}
		\caption{Memory usage.}
		\label{fig:hel3}
	\end{subfigure}
	\begin{subfigure}{0.4\textwidth}
		\centering
		\includegraphics[width=\linewidth]{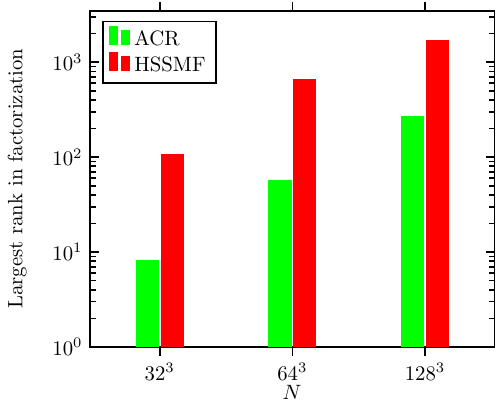}
		\caption{Largest rank in the factorization.}
		\label{fig:hel4}
	\end{subfigure}
\caption{Solution of increasingly larger indefinite Helmholtz problems consistently discretized with 12 points per wavelength.}
\label{fig:helm}
\end{figure}

We finally consider the indefinite Helmholtz equation with Dirichlet boundary conditions on the unit cube, i.e.

\begin{equation}
\begin{aligned}
- (\nabla^2u + \kappa^2u) = 1, \; \Omega = [0,1]^3,  \\
\end{aligned}
\label{helmholtzEquation}
\end{equation}

\noindent discretized with the 27-point trilinear finite element scheme on hexahedra. Results for ACR and HSSMF are reported in Figure \ref{fig:helm}. The parameter $\kappa$ is chosen to obtain a sampling rate of approximately 12 points per wavelength, specifically $\kappa=\{16,32,64\}$ respectively, corresponding to approximately $10 \times 10 \times 10$ for the $N=128^3$ problem. As opposed to the positive definite Helmholtz equation which models phenomena similar to diffusion, the indefinite variant, commonly denoted as the wave Helmholtz equation, has a solution that is oscillatory in nature. Multigrid methods are known to diverge without specific customizations for high-frequency Helmholtz problems, which we also confirmed via experimentation. For a detailed examination of the difficulties of solving the Helmholtz equation with classical iterative methods we refer the reader to \cite{ernst2012difficult}.

We document the execution parameters, obtained relative residual, and ranks of the ACR and HSSMF factorization in Tables \ref{table:ACRHelConfiguration} and \ref{table:STRHelmConfiguration}. Numerical experiments show that ACR features consistently lower factorization and solve times than HSSMF, as can be seen in Figure \ref{fig:hel1} and \ref{fig:hel2}. The size of the factors of ACR and HSSMF are comparable, with a slightly higher memory requirements of ACR due to performance-oriented tuning, see Figure \ref{fig:hel3}. Furthermore, as also shown in section \ref{section:poisson}, HSSMF required less memory than MF, and CR quickly runs out of memory for problems larger than $N=64^3$. Finally, the largest rank of ACR is consistently lower than that of HSSMF, even though both solvers require $O(n)$ ranks, as shown in Figure \ref{fig:hel4}. Nevertheless, lower ranks lead to faster time-to-solution in favor of ACR.

\begin{table}[]
\centering
\begin{tabular}{|c|c|l|c|c|c|c|}
\hline
\textbf{$N$} & \textbf{$\mathcal{H}_{\epsilon}$} & \multicolumn{1}{c|}{\textbf{$\eta$}} & \textbf{\begin{tabular}[c]{@{}c@{}}Leaf\\ size\end{tabular}} & \textbf{\begin{tabular}[c]{@{}c@{}}Relative\\ residual\end{tabular}} & \multicolumn{1}{l|}{\textbf{Average rank}} & \multicolumn{1}{l|}{\textbf{Largest rank}} \\ \hline
$32^3$       & 5E-03                          & 4                                    & 32                                                           & 1.67E-02                                                             & 5                                          & 8                                          \\ \hline
$64^3$       & 5E-08                          & 8                                    & 32                                                           & 2.63E-02                                                             & 30                                         & 56                                         \\ \hline
$128^3$      & 5E-13                          & 16                                   & 32                                                           & 1.07E-02                                                             & 113                                        & 260                                        \\ \hline
\end{tabular}
\caption{Execution parameters, obtained relative residual, and ranks of the ACR factorization for the Helmholtz experiments.}
\label{table:ACRHelConfiguration}

\bigskip

\centering
\begin{tabular}{|c|c|c|c|c|c|c|}
\hline
\textbf{$N$} & \textbf{\begin{tabular}[c]{@{}c@{}}Compression\\ tolerance\end{tabular}} & \multicolumn{1}{l|}{\textbf{\begin{tabular}[c]{@{}l@{}}Relative\\ residual\end{tabular}}} & \textbf{\begin{tabular}[c]{@{}c@{}}Leaf\\ size\end{tabular}} & \textbf{\begin{tabular}[c]{@{}c@{}}Minimum\\ front size\end{tabular}} & \multicolumn{1}{l|}{\textbf{\begin{tabular}[c]{@{}l@{}}Largest rank\end{tabular}}} \\ \hline
$32^3$       & 5E-03                                                                    & 5.32E-02                                                                                  & 128                                                          & 256                                                                   & 105                                                                                                   \\ \hline
$64^3$       & 1E-04                                                                    & 6.08E-02                                                                                  & 128                                                          & 1,024                                                                  & 641                                                                                                   \\ \hline
$128^3$      & 1E-06                                                                    & 1.13E-02                                                                                  & 128                                                          & 4,096                                                                  & 1,659                                                                                                  \\ \hline
\end{tabular}
\caption{Execution parameters, obtained relative residual, and ranks of the HSSMF factorization for the Helmholtz experiments.}
\label{table:STRHelmConfiguration}
\end{table}

\clearpage
\section{Conclusions and future work}

We present a novel fast direct solver, Accelerated Cyclic Reduction, for block tridiagonal linear systems which commonly arise in the discretization of elliptic operators. The elimination strategy is based on a red/black ordering of the blocks that logically divides the grid into planes, approximates matrix blocks representing these planes with $\mathcal{H}$-matrices, and proceeds with elimination using hierarchical matrix operations. ACR achieves log-linear arithmetic complexity of $O(k~N \log N~(\log N + k^2))$ and memory requirements of $\mathcal{O}(k~N \log N)$ by approximating each block with a hierarchical matrix whose structure is refined using a spatial partitioning of the planar grid sections, employing a strong admissibility criterion that effectively limits the ranks of individual low rank blocks in the hierarchical matrix representations, and operating with hierarchical matrix arithmetics throughout. The average rank $k$ of the blocks inside the hierarchical matrix representations controls the accuracy of the approximation and grows only modestly with problem size. A fair agreement with the rank estimate of \cite{chandrasekaran2010} was found for the 3D Poisson equation of $O(1)$ (Table \ref{table:ACRPoisConfiguration}), and for the 3D Helmholtz equation $O(n)$ (Table \ref{table:ACRHelConfiguration}).

The concurrency features of ACR are among its most important strengths. The regularity and structure of the decompositions allow efficient load balance. These features are demonstrated in a distributed-memory environment with numerical experiments that study the strong and weak scalability of our implementation. We provide a reference for performance and memory consumption using comparisons with state-of-the-art open-source implementations of the HSS-structured multifrontal solver from the STRUMPACK library, and algebraic multigrid from hypre. 

ACR, being essentially a direct solver with tunable accuracy, can tackle problems that lack definiteness, such as the indefinite high-frequency Helmholtz equation, or symmetry, such as the convection-diffusion equation.
For these problems, stock versions of algebraic multigrid fail to produce convergent schemes. We demonstrated the robustness of ACR in dealing with such problems over a range of problem sizes and parameters.

While multigrid methods are generally superior for scalar problems possessing smoothness and definiteness, direct factorization methods such as ACR and HSSMF benefit where multiple right-hand sides are involved, as the time to solve per extra forcing term is orders of magnitude smaller than the factorization, which can be reused. The smaller ranks $k$ of ACR result in solution times per new right-hand side that are smaller than those of HSSMF.

Although having the same asymptotic complexity as other solvers that use general hierarchical matrix representations in their factorizations, such as $\mathcal{H}$-LU, ACR has fundamentally different algorithmic roots which enable a novel alternative for a relevant class of problems with competitive performance, increasing concurrency as the problem grows and almost optimal memory requirements. Moreover, to the best of our knowledge, this is the first distributed-memory implementation of the synergies of cyclic reduction and hierarchical matrices, which scales up to 8,192 cores for problems up to $N=512^3$ degrees of freedom.

\begin{figure}[H]
	\centering
	\begin{subfigure}{0.4\textwidth}
		\centering
		\includegraphics[width=0.7\linewidth]{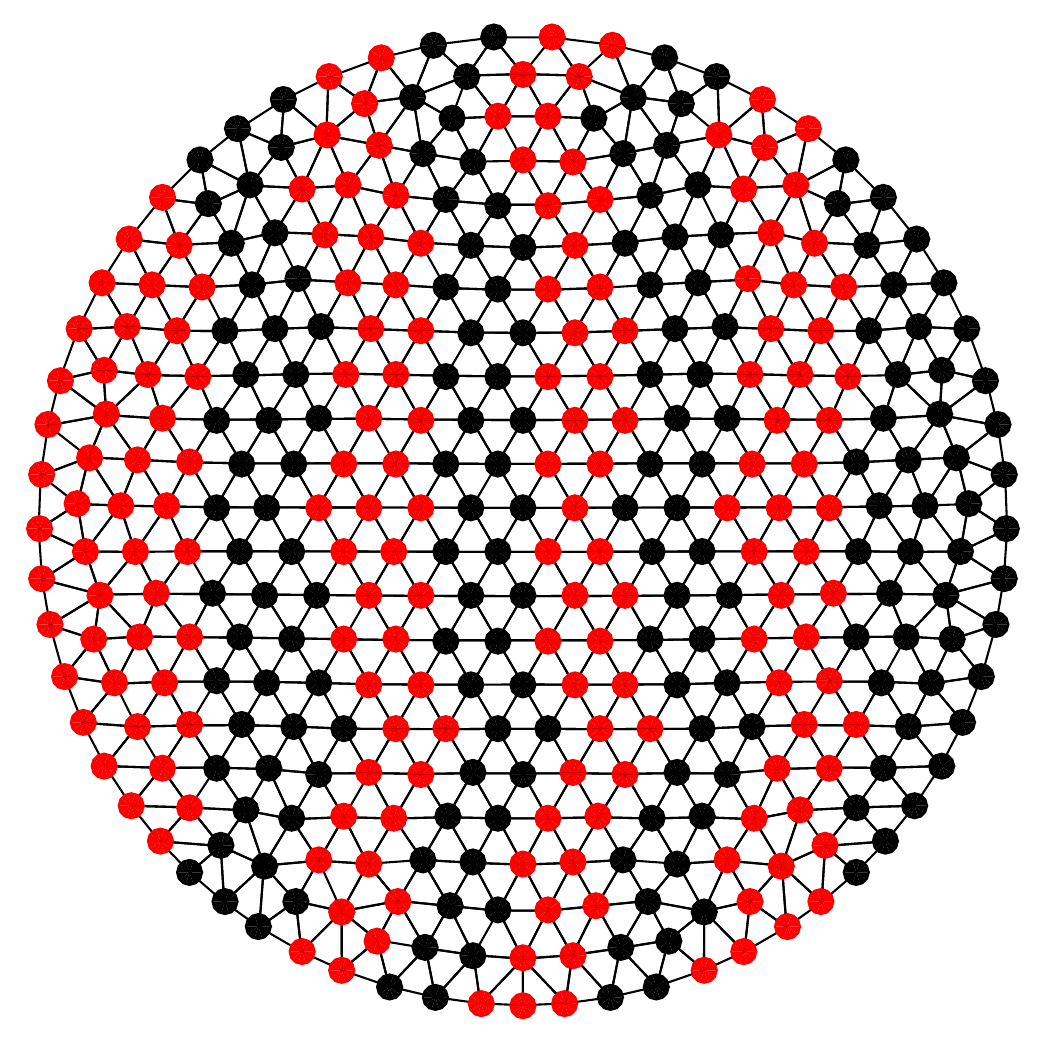}
	\end{subfigure}
	\begin{subfigure}{0.4\textwidth}
		\centering
		\includegraphics[width=0.7\linewidth]{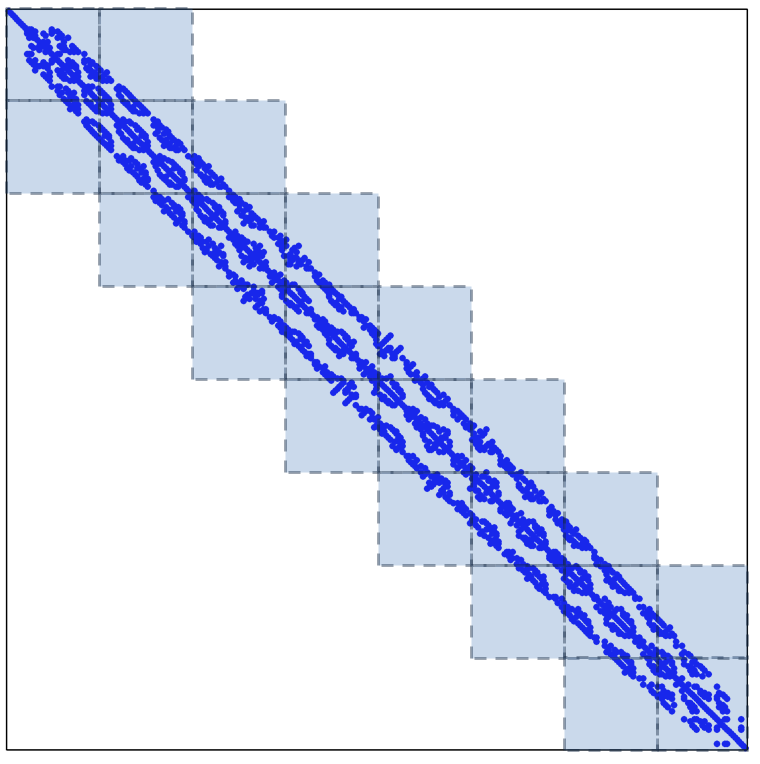}
	\end{subfigure}
\caption{Partitioning of an unstructured mesh that produces a block tridiagonal matrix structure, for the application of ACR.}
\label{fig:unsACR}
\end{figure}

ACR has been demonstrated for a regular grid discretization, but the generalization to arbitrary grids is possible and we intend to explore it in the future. Such a generalization would require an ordering of the mesh that produces a sequence of thin elongated regions (in 2D or 3D) where every region has only two neighbors so that the block tridiagonal structure is preserved. Such an ordering might be produced via a breadth-first search traversal of the mesh as shown in Figure \ref{fig:unsACR}. In the unstructured case, the diagonal blocks do not necessarily have the same size, and the off-diagonal blocks might be of rectangular shape. The main algorithmic implication is that each block will now have its own hierarchical matrix structure generated from the geometry of the region it represents. Computationally however, the structure generation represents a small portion in the overall computation.

In addition, because of the tunable accuracy characteristics of ACR, there are complexity-accuracy trade-offs that would naturally lead to the development of a new scalable preconditioner which we present at \cite{Chavez2017}.

\section{Acknowledgments}
We thank the anonymous reviewers for their detailed comments and suggestions for this manuscript. The authors would also like to thank Ronald Kriemann from the Max-Planck-Institute for Mathematics in the Sciences for development and continuous support of HLibPro, Alexander Litvinenko from the King Abdullah University of Science and Technology (KAUST) for the enlightening discussions and advice, and Pieter Ghysels from the Lawrence Berkeley National Laboratory for his recommendations on the use of STRUMPACK. Support from the KAUST Supercomputing Laboratory and access to Shaheen is gratefully acknowledged. The work of all authors was supported by the Extreme Computing Research Center at KAUST.

\bibliographystyle{elsarticle-num}
\bibliography{PARCO_DirectACR}

\vspace{1in}
\textbf{BibTeX entry of this article:}
\begin{verbatim}
@article{Chavez2016,
author = "Gustavo Ch{\'a}vez and George Turkiyyah and Stefano Zampini
and Hatem Ltaief and David Keyes",
title = "Accelerated Cyclic Reduction: A Distributed-Memory Fast Solver
for Structured Linear Systems",
journal = "Parallel Computing",
year = "2017",
issn = "0167-8191",
doi = "https://doi.org/10.1016/j.parco.2017.12.001",
url = "https://www.sciencedirect.com/science/article/pii/S0167819117302041",
}
\end{verbatim}

\end{document}